\documentclass[12pt]{article}
\usepackage{amssymb,amsmath,graphicx,acronym,amsfonts,float,color,epstopdf,longtable,subfigure,caption2,makecell,amsthm,bm}
\usepackage{amsfonts,ifpdf,dsfont,longtable}
\usepackage[noend]{algpseudocode}
\usepackage{amsmath}
\usepackage{algorithm}
\usepackage{algcompatible}
\usepackage{setspace}
\DeclareMathOperator{\Tr}{Tr}

\newtheorem{theorem}{Theorem}[section]

\newtheorem{definition}{Definition}

\newtheorem{Remark}{Remark}[section]

\def\proof{{\noindent\bf Proof.~~}}

\def\R{\mathbb{R}}
\def\A{\textbf{{A}}}

\def\H{\textbf{{H}}}
\def\K{\textbf{{K}}}
\def\M{\textbf{{M}}}
\def\G{\textbf{{G}}}
\def\U{\textbf{{U}}}
\def\B{\textbf{{B}}}
\def\X{\textbf{{X}}}
\def\Y{\textbf{{Y}}}

\def\e{{\bf e}}

\def\x{\textbf{\textit{x}}}
\def\y{\textbf{\textit{y}}}
\def\w{\textbf{\textit{w}}}

\def\u{\textbf{\textit{u}}}

\def\0{{\bf 0}}

\def\be{\begin{equation}}
\def\ee{\end{equation}}
\input epsf
\date{}  
\begin{document}

\title{\bf A Nonparallel Support Tensor Machine for Binary Classification based Large Margin Distribution and Iterative Optimization\thanks{The first author's work was supported by the National Natural Science Foundation of P.R. China (Grant No.12171064), by The team project of innovation leading talent in chongqing (No.CQYC20210309536) and by the Foundation of Chongqing Normal university (20XLB009)}.
}

\author{Zhuolin Du, Yisheng Song\thanks{Corresponding author E-mail: yisheng.song@cqnu.edu.cn}\\
School of Mathematical Sciences,  Chongqing Normal University, \\
Chongqing, 401331, P.R. China. \\ Email: duzhuolin728@163.com (Du); yisheng.song@cqnu.edu.cn (Song)}

\maketitle

\begin{abstract}

Based on the tensor-based large margin distribution and the nonparallel support tensor machine, we establish a novel classifier for binary classification problem in this paper, termed the Large Margin Distribution based NonParallel Support Tensor Machine (LDM-NPSTM).  The proposed classifier has the following advantages: First, it utilizes tensor data as training samples, which helps to comprehensively preserve the inherent structural information of high-dimensional data, thereby improving classification accuracy. Second, this classifier not only considers traditional empirical risk and structural risk but also incorporates the marginal distribution information of the samples, further enhancing its classification performance. To solve this classifier, we use alternative projection algorithm. Specifically, building on the formulation where in the proposed LDM-NPSTM, the parameters defining the separating hyperplane form a tensor (tensorplane) constrained to be the sum of rank-one tensors, the corresponding optimization problem is solved iteratively using alternative projection algorithm. In each iteration, the parameters related to the projections along a single tensor mode are estimated by solving a typical Support Vector Machine-type optimization problem. Finally, the efficiency and performance of the proposed model and algorithm are verified through theoretical analysis and some numerical examples.\\



\end{abstract}
\newpage

\section{Introduction}
A significant number of real-world datasets, especially those involving image data, are frequently represented in tensor format. For instance, a grayscale face image \cite{kr97} can be modeled as a second-order tensor (or matrix), while color images \cite{gd99,pk01}, grayscale video sequences \cite{rl07}, gait contour sequences \cite{hk08}, and hyper spectral cubes \cite{nb09} are typically expressed as third-order tensors. Additionally, color video sequences are often represented as fourth-order tensors \cite{km01,wa04}. The tensor-based data representation, with its multi-dimensional structure that complicates the capture of spatial and temporal relationships, high dimensionality prone to the curse of dimensionality and overfitting, and substantial storage and computational requirements, uniquely complicates feature extraction and representation, thereby rendering it a fundamental challenge in the design of classifier models.

One of the most representative and successful classification algorithms is the Support Vector Machines (SVM) \cite{cvv95,ss96,vv13}, which have been successfully applied to a variety of real-world pattern recognition problems, such as text classification \cite{il08,an10}, image classification \cite{sj10,ra11}, feature extraction \cite{lc08,sk08,ds08}, web mining \cite{by10} and function estimation \cite{cf03,nl10}. The central idea of SVM is to find the optimal separating hyperplane between the positive and negative examples. The optimal hyperplane is defined as the one giving maximum margin between the training examples that are closest to the hyperplane. Different from the traditional SVM, in 2007, Jayadeva et al. \cite{krs07} proposed Twin SVM (TWSVM), which also aims at generating two nonparallel planes such that each plane is closer to one of the two classes and is as far as possible from the other. Notably, the formulation of TWSVMs is very much in line with standard SVMs. However, TWSVMs seek to solve two dual QPPs of smaller size rather than solving single dual QPP with large number of parameters in conventional SVM. As a result, the algorithm achieves a processing speed roughly fourfold faster compared to the traditional SVM \cite{krs07}. While the aforementioned classification methods focus on maximizing the minimum margin, researchs by Gao et al. \cite{gz13,lm11} have indicated that doing so does not necessarily guarantee improved generalization performance. Instead, the distribution of margins has been shown to play a more critical role. Here, the margin distribution is defined as the margin mean and the margin variance. Therefore, in order to improve the generalization performance of SVM, Zhou et al. \cite{zz14} characterized margin distribution through its mean and variance, leading to the development of the Large Margin Distribution Machine (LDM), which builds upon the SVM framework. The effectiveness of LDM has been proved in theory and experiments.

In recent years, there has been a growing interest in extending traditional vector or matrix-based machine learning algorithms to better handle tensor data \cite{cc22,hz20}. This shift is motivated by the need to effectively process high-dimensional datasets, such as those encountered in image and video analysis. In 2005, Tao et al. \cite{dx05} proposed a Supervised Tensor Learning (STL) scheme by replacing the vector inputs with tensor inputs and decomposing the corresponding weight vector into a rank-1 tensor, which is trained by the alternating projection optimization method. Based on this learning scheme, in 2007, Tao et al. \cite{dx07} further extended the standard linear SVM to a tensorial format known as the Support Tensor Machine (STM). This adaptation allows for more effective classification of tensor data by leveraging its inherent structure. Following this development, Zhang et al. \cite{zgw09} generalized the vector-based learning algorithm TWSVM to the tensor-based method Twin STM (TWSTM), and implemented the classifier for microcalcification clusters detection. By comparison with TWSVM, the tensor version reduces the overfitting problem significantly. Additionally, Khemchandani et al. developed a least squares variant of STM, termed as Proximal STM (PSTM) \cite{kk13}, where the classifier is obtained by solving a system of linear equations rather than a quadratic programming problem at each iteration of PSTM algorithm as compared to STM algorithm.  This modification enhances computational efficiency while maintaining classification performance. Tensor-based algorithms on the other hand decompose the whole problem into several smaller and simpler subproblems, each defined over specific tensor modes and characterized by lower dimensionality. This decomposition has been shown to reduce the degree of overfitting that appears in vector-based learning techniques, particularly when few training samples are available \cite{dx07}.

In this paper, we propose a novel framework, termed Nonparallel Support Tensor Machine based on Large Margin Distribution (LDM-NPSTM), aimed at further enhancing the generalization performance of the Twin Support Tensor Machine (TWSTM). Drawing on the strengths of Large Margin Distribution theory \cite{zz14} and TWSTM \cite{zgw09}, our approach integrates their core principles to address classification challenges more effectively. Specifically, we characterize the margin distribution using first-order (margin mean) and second-order (margin variance) statistics, with the core objective of maximizing the margin mean while minimizing the margin variance to improve classification robustness.
To ensure a more rigorous model structure, we incorporate a regularization term into the LDM-NPSTM framework, balancing empirical risk and structural complexity. For model optimization, we adopt an iterative solution based on CANDECOMP/PARAFAC (CP) decomposition. In each iteration, parameters corresponding to projections along a single tensor mode are estimated by solving a typical SVM-type optimization problem. Notably, the inverse matrix involved in the dual problem is inherently nonsingular, eliminating the need for additional assumptions and simplifying the computational process.

The remainder of this paper is structured as follows: In Section 2, we introduce the notations consistently used throughout the paper and provide a concise overview of fundamental concepts, including those related to SVM, TWSVM, TWSTM and LDM. Section 3 elaborates on our proposed framework, the LDM-NPSTM, with detailed formulations and a discussion of its key advantages. Experimental results that demonstrate the effectiveness of the LMD-NPSTM are discussed in Section 4. At last, concluding remarks are given in Section 5.

\section{Preliminaries}
In this section, we first introduce some notation and basic definitions used throughout the paper, and then briefly review the related works.
\subsection{Notation and Basic Definitions}
An \( m \)-th order tensor is defined as a collection of measurements indexed by \( m \) indices, with each index corresponding to a mode. Vectors are considered first-order tensors, while matrices represent second-order tensors \cite{kt09}. In this paper, we will utilize lowercase letters (e.g., \( x \)) to denote scalars, boldface lowercase letters (e.g., \(\mathbf{x}\)) and boldface capital letters (e.g., \(\mathbf{X}\)) to represent vectors and matrices, respectively. Tensors of order 3 or higher will be denoted by boldface Euler script calligraphic letters (e.g., \(\mathcal{X}\)). Furthermore, we denote the set of all mth-order n-dimensional real tensors as $\mathrm{T}_{m,n}$. The \( i \)-th element of a vector \( \mathbf{x} \in \mathbb{R}^n \) is denoted by \( x_i, \, i = 1, 2, \ldots, n\). In a similar way, the elements of an \( m \)-th order tensor \( \mathcal{X} \) will be denoted by \( x_{i_1 i_2 \ldots i_m}\), where \(i_j= 1, 2, \ldots,n_j\) for \(j=1,\ldots,m\). Moreover, we summarize some notations used throughout the paper in Table 1.

\begin{table}[ht]
    \centering
    \caption{List of Symbols}
    \begin{tabular}{@{}ll@{}}
        \hline
        Symbol & Description \\ \hline
        $m$ & The total number of samples \\
        $m_1$ & The number of positive samples \\
        $m_2$ & The number of negative samples \\
        $\mathcal{X}_p$ & The $p$th input positive tensor sample \\
        $\mathcal{Y}_q$ & The $q$th input negative tensor sample \\
        $y_i$ & The label of training samples\\
        $M$  & The order of $\mathcal{X}_p, \mathcal{Y}_q\in \mathbb{R}^{I_1 \times I_2 \times \cdots \times I_M}$ \\
        $x_i$ & The $i$th input vector sample \\
        $\mathbf{w}, \mathcal{W}$ & The weight parameters of vector and tensor \\
        $\alpha, \beta$ & The Lagrange multipliers \\
        $\mathbf{w}^{(1)} \circ \mathbf{w}^{(2)} \circ \cdots \circ \mathbf{w}^{(M)}$ & Rank-one tensor of $\mathcal{W}$ \\
        $\mathcal{X} \times_n \mathbf{w}$ & n-mode product of $\mathcal{X}$ and $\mathbf{w}$ \\
        $\|\cdot\|_F$ & Frobenius norm \\ \hline
    \end{tabular}
    \label{tab:symbols}
\end{table}

In the following, we introduce some notation and definitions of the tensors and matrices in the area of multilinear algebra \cite{dl97,kt09}.

\begin{definition}
\textnormal{(Inner product)} Given tensors \( \mathcal{X}, \mathcal{Y} \in \mathbb{R}^{I_1 \times \cdots \times I_M} \), the inner product of \( \mathcal{X} \) and \( \mathcal{Y} \) is defined as
\begin{equation}
\langle\mathcal{X},\mathcal{Y}\rangle := \sum_{i_1=1}^{I_1} \sum_{i_2=1}^{I_2} \cdots \sum_{i_M=1}^{I_M} x_{i_1 i_2 \ldots i_M} y_{i_1 i_2 \ldots i_M}.
\end{equation}\label{2.1}
\end{definition}
\begin{definition}
\textnormal{(Frobenius norm)} The Frobenius norm of a tensor \( \mathcal{A} \in \mathbb{R}^{I_1 \times I_2 \times \cdots \times I_M} \) is defined as
\begin{equation}\label{2.2}
\|\mathcal{A}\|_F := \sqrt{\langle \mathcal{A}, \mathcal{A} \rangle}.
\end{equation}
\end{definition}

\begin{Remark}
Given two same-sized tensors \( \mathcal{A} \in \mathbb{R}^{I_1 \times I_2 \times \cdots \times I_M} \) and \( \mathcal{B} \in \mathbb{R}^{I_1 \times I_2 \times \cdots \times I_M} \), the distance between tensors \( \mathcal{A} \) and \( \mathcal{B} \) is defined as \( \| \mathcal{A}-\mathcal{B} \|_F \). Note that the Frobenius norm of the difference between two tensors equals the Euclidean distance of their vectorized representations \cite{lp09}.
\end{Remark}

\begin{definition}\label{2.3}
\textnormal{(Outer product)} We use $\otimes$ to denote tensor outer product; that is, for any two tensors $\mathcal{A}\in \mathrm{T}_{m,n}$ and $\mathcal{B}\in\mathrm{T}_{p,n}$ is given by:
\begin{equation}\label{2.3}
\mathcal{A}\otimes\mathcal{B}=(a_{i_1\ldots i_m}b_{i_1\ldots i_p})\in\mathrm{T}_{m+p,n}.
\end{equation}
\end{definition}
According to this definition, it is easy to check that
\[\x^{{\otimes}k}=\underbrace{\x\otimes\cdots\otimes\x}_{k~ \mbox{times}}=(x_{i_1}\cdots x_{i_k})\in\mathrm{T}_{k,n}.\]

\begin{definition}\label{2.4}
\textnormal{(CP decomposition)} Given \( \mathcal{X} \in \mathbb{R}^{I_1 \times \cdots \times I_M} \), if there exist \( \u_r^{(1)} \in \mathbb{R}^{I_1} \), \( \u_r^{(2)} \in \mathbb{R}^{I_2} \), \ldots, \( \u_r^{(M)} \in \mathbb{R}^{I_M} \) such that
\begin{equation}\label{2.4}
\mathcal{X} = \sum_{r=1}^{R} \u_r^{(1)} \otimes \u_r^{(2)} \otimes \cdots \otimes \u_r^{(M)},
\end{equation}
where R is a positive integer, we call (\ref{2.4}) a tensor CANDECOMP/PARAFAC (CP) decomposition of \( \mathcal{X} \).
\end{definition}

\begin{definition}\label{3.5}
\textnormal{(Matricization)}
The matricization (also known as unfolding or flattening of a tensor) is the reordering of the tensor elements into a matrix. The \( n \)-mode matricization of a tensor \( \mathcal{X} \in \mathbb{R}^{I_1 \times I_2 \times \cdots \times I_M} \), denoted by \( \mathcal{X}_{(n)} \in \mathbb{R}^{I_n \times \left(\prod_{k \neq n} I_k\right)} \), arranges the \( n \)-mode fibers to become the columns of the final matrix. Each tensor element \( (i_1, i_2, \ldots, i_M) \) maps to the matrix element \( (i_n, j) \), where

\[
j = 1 + \sum_{k=1, k \neq n}^{M} (i_k - 1) J_k \quad \textnormal {with} \quad J_k = \prod_{l=1, l \neq n}^{k-1} I_l.
\]
\end{definition}

A more general treatment of matricization can be found in \cite{kt06}.

\begin{definition}\label{3.6}
\textnormal{(Matrix kronecker product)}
The Kronecker product of matrices \( \A \in \mathbb{R}^{I \times J} \) and \( \B \in \mathbb{R}^{K \times L} \) is denoted by \( \A \otimes \B \). The result is a matrix of size \( (IK) \times (JL) \) and defined by
\[
\A \otimes \B =
\begin{bmatrix}
a_{11}\B & a_{12}\B & \cdots & a_{1J}\B \\
a_{21}\B & a_{22}\B & \cdots & a_{2J}\B \\
\vdots & \vdots & \ddots & \vdots \\
a_{I1}\B & a_{I2}\B & \cdots & a_{IJ}\B
\end{bmatrix}
\]
\end{definition}

\begin{definition}\label{def3.7}
\textnormal{(Matrix Khatri-Rao product)}
Given matrices \( \A \in \mathbb{R}^{I \times K} \) and \( \B \in \mathbb{R}^{J \times K} \), their Khatri-Rao product is denoted by \( \A \odot \B \). The result is a matrix of size \( (IJ) \times K \) defined by
\[
\A \odot \B = [\bm a_1 \otimes \bm b_1 \quad \bm a_2 \otimes \bm b_2 \quad \cdots \quad \bm a_K \otimes \bm b_K ].
\]

\end{definition}

\begin{Remark}
If matrices \( \A \) and \( \B \) of Definition \ref{def3.7} are vectors, i.e., $\bm a$ and $\bm b$, then the Khatri-Rao and Kronecker products are identical, i.e., \( \bm a\otimes\bm b=\bm a\odot\bm b \).
\end{Remark}

\subsection{Related works}
Support Vector Machine (SVMs) form a class of supervised machine learning algorithms that train the classifier function using pre-labeled data. Specifically, for a given training set $\{(\x_i, y_i)\mid i=1,\cdots,m\}$, where data points $\x_i\in \R^n$ and the class label $y_i=\{-1,1\}$, the objective of the support vector machine problem is to identify a hyperplane $\w^\top \x+b=0$, where $\w \in \R^n$ and $b\in \R$, in such a way that the two different classes of data points are separated with a maximal separation margin and minimal the classification loss. The standard SVM problem can be formulated as the following convex quadratic problem \cite{vv99}:
\begin{equation}\label{2.5}
\begin{aligned}
\min_{\w,b,\xi_i}\quad&\frac12\|\w\|^2+C\sum_{i=1}^{m}\xi_i\\
{\rm s.t.}\quad&y_i(\w^\top \x_i+b)\geq 1-\xi_i,\quad i=1,2,\ldots,m\\
&\xi_i\geq0, \quad i=1,2,\ldots,m.
\end{aligned}
\end{equation}
where $\xi_i$ is a slack variable, and $C>0$ is a penalty parameter that represents the loss weight.

Moreover, the standard SVM model requires solving a single large-scale optimization problem, which can be computationally intensive. To address this issue, Khem-chandani et al. proposed Twin SVM (TWSVM) \cite{krs07}. TWSVM generates two non-parallel planes such that each plane is closer to one of the two cases and is as far as possible from the other. This approach allows TWSVM to solve a pair of smaller-sized quadratic programming problems (QPPs) rather than a single large QPP, resulting in a computational speed that is approximately four times faster than that of traditional SVMs. The optimization problems in TWSVM can be formulated as the following pair of quadratic programming problems:
\begin{equation}\label{2.6}
\begin{aligned}
\min_{\w^{(1)}, b^{(1)}, \bm q} \quad&\frac{1}{2} (\A\w^{(1)}+\e_1b^{(1)})^{\top}(\A\w^{(1)}+\e_1b^{(1)})+ c_1\e_2^\top \bm q \\
\rm{s.t.}\quad& -(\B\w^{(1)} + \e_2b^{(1)})+\bm q\geq \e_2,\quad\bm q \geq \0\\
\end{aligned}
\end{equation}
and
\begin{equation}\label{2.7}
\begin{aligned}
\min_{\w^{(2)}, b^{(2)}, \bm q} \quad&\frac{1}{2}(\B\w^{(2)}+ \e_2b^{(2)})^{\top}(\B\w^{(2)}+ \e_2 b^{(2)})+c_2\e_1^\top\bm q  \\
\rm{s.t.}\quad& -(\A\w^{(2)}+\e_1b^{(2)})+\bm q \geq\e_1,\quad \bm q \geq \0\\
\end{aligned}
\end{equation}
where $c_1,c_2>0$ are penalty parameters, $\e_1$ and $\e_2$ are the vectors of one of appropriate dimensions, $\bm q\in\mathbb{R}^n$ is the vector of all slack variables to deal with linearly nonseparable problems.

Both SVM and TWSVM are vector-based learning algorithms, which accept vectors as inputs. In practice, real-world image and video data are more naturally represented as matrices (second-order tensors) or higher-order tensors. Therefore, Zhang et al. \cite{zgw09} generalize the vector-based learning algorithm TWSVM to the tensor-based method Twin Support Tensor Machines (TWSTM), which accepts general tensors as input. The following formulation for TWSTM can be established
\begin{equation}\label{2.8}
\begin{aligned}
\min_{\w^{(1)}_k, b^{(1)}, \bm q}\quad&\frac{1}{2} \left( \mathcal{A} \prod_{k=1}^{K} \times_k \w^{(1)}_k + \e_1 b^{(1)} \right)^{\top} \left( \mathcal{A} \prod_{k=1}^{K} \times_k \w^{(1)}_k + \e_1 b^{(1)} \right) + c_1 \e_2^{\top} \bm q \\
\rm{s.t.} \quad&-\left( \mathcal{B} \prod_{k=1}^{K} \times_k \w^{(1)}_k + \e_2 b^{(1)} \right) + \bm q \geq \e_2, \quad \bm q \geq \0
\end{aligned}
\end{equation}
and
\begin{equation}\label{2.9}
\begin{aligned}
\min_{\w^{(2)}_k, b^{(2)},\bm q} \quad& \frac{1}{2} \left( \mathcal{B} \prod_{k=1}^{K}\times_k \w^{(2)}_k + \e_2 b^{(2)} \right)^{\top} \left( \mathcal{B} \prod_{k=1}^{K} \times_k \w^{(2)}_k + \e_2 b^{(2)} \right) + c_2 \e_1^{\top} \bm q \\
\rm{s.t.} \quad& \mathcal{A} \prod_{k=1}^{K} \times_k \w^{(2)}_k + \e_1 b^{(2)} + \bm q \geq \e_1, \quad \bm q \geq \0
\end{aligned}
\end{equation}
where $c_1$ and $c_2$ are penalty parameters, $\e_1$ and $\e_2$ are the vectors of appropriate dimensions, $\bm q\in\mathbb{R}^n$ is the vector of all slack variables to deal with linearly nonseparateable problems.

The objective function of the above model is based on minimizing the margin. From the perspective of structural risk, Gao et al. \cite{gz13} have verified that the margin distribution is more important than minimum margins in optimizing generalization performance. By characterizing the margin distribution in terms of margin mean and margin variance, Zhou et al. \cite{zz14} proposed the Large Margin Distribution Machine (LDM) on the basis of SVM, which focuses on optimizing the distances from the center of the other category. The following formulation for LDM can be established
\begin{equation}\label{2.10}
\begin{aligned}
\min_{\w, \bm \xi} &\quad \frac{1}{2} \w^{\top}\w + \lambda_1 \hat{\gamma}-\lambda_2\bar{\gamma}+C \sum_{i=1}^{m} \bm\xi_i \\
\rm{s.t.} & \quad y_i \w^{\top} \x_i \geq 1 - \xi_i, \\
& \quad \xi_i \geq 0, \quad i = 1,\ldots, m.
\end{aligned}
\end{equation}
where $\lambda_1$ and $\lambda_2$ are the parameters for trading-off the margin variance $\hat{\gamma}$, the margin mean $\bar{\gamma}$ and the model complexity.

%


\setcounter{proposition}{0} \setcounter{theorem}{0}
\setcounter{lemma}{0} \setcounter{corollary}{0}
\setcounter{equation}{0}
\setcounter{definition}{0}

\section{Proposed LMD-NPSTM }
In this section, LDM-NPSTM is proposed. Specifically, Subsection 3.1 introduces the structure of the model, optimization strategies, and associated dual problems. Moreover, the detailed implementation algorithm of LMD-NPSTM will be shown in Subsection 3.2.

\subsection{Model construction and optimization}
Consider the binary classification problem in tensor space, where the training set is defined as
$\mathcal{T}_m = \{(\mathcal{X}_p, y_p)|p=1,2,\ldots,m_1\}\cup\{ (\mathcal{Y}_q, y_q)|q =1,2,\ldots,m_2\}$. Here, $\mathcal{X}_p, \mathcal{Y}_q \in\mathbb{R}^{I_1\times I_2\times\cdots\times I_M}$ denotes the feature tensor of the $p,q$-th sample, respectively, $y_p=1$ and $y_q=-1$ are their corresponding class labels; and $m = m_1 + m_2$ (where $m_1$ and $m_2$ are the numbers of positive and negative samples in $\mathcal{T}_m$, respectively).
Let $\y_1 = [1, \cdots, 1]^\top \in \mathbb{R}^{m_1 \times 1}$ be the label vector for all positive samples, where each element corresponds to the label of a positive instance, and $\mathbf{y}_2 = [-1, \cdots, -1]^\top \in \mathbb{R}^{m_2 \times 1}$ denote the label vector for all negative samples, with each element corresponding to the label of a negative instance. The LMD-NPSTM identifies two non-parallel hyperplanes in the feature space:


$$f_1(\mathcal{X})=
\langle\mathcal{W}_{1},\mathcal{X}\rangle=0$$
$$f_2(\mathcal{Y})=
\langle\mathcal{W}_{2},\mathcal{Y}\rangle=0,$$
where $\mathcal{W}_{1},\mathcal{W}_{2}\in \mathbb{R}^{I_1\times I_2\times\cdots\times I_M}$.

Notably, inspired by the work of Zhou et al. \cite{zz14}, the bias term does not affect the overall derivation process.
Furthermore, unlike the conventional approach that calculates margin mean and variance using the entire dataset, the LMD-NPSTM separates these calculations into positive-class and negative-class components \cite{zl23}. Specifically, these metrics are associated with the distances to their respective hyperplanes, as elaborated below.

The distance from an individual data point to the hyperplane is defined by
\begin{equation*}
\begin{cases}
\gamma_+^p=y_{1}^p f_2\big(\mathcal{X}_{p}\big)=y_{1}^p\langle\mathcal{W}_{2}, \mathcal{X}_{p}\rangle,~~ p=1,\ldots, m_1\\
\gamma_-^q=y_{2}^q f_1\big(\mathcal{Y}_{q}\big)=y_{2}^q\langle\mathcal{W}_{1}, \mathcal{Y}_{q}\rangle,~~ q=1,\ldots, m_2,
\end{cases}
\end{equation*}
where $y_1^p \in \{1, -1\}$ and $y_2^q \in \{1, -1\}$ denote the labels of positive/negative samples, respectively.\\
The margin means are defined as follows:

\begin{equation}\label{3.1}
\begin{cases}
\bar{\gamma}_{+}=\frac{1}{m_{1}}\sum\limits_{p=1}^{m_{1}}\gamma_{+}^{p}=\frac{1}{m_{1}}\sum\limits_{p=1}^{m_{1}} y_{1}^{p}f_{2}(\mathcal{X}_{p}) =\frac{1}{m_{1}}y_{1}^{p}\langle \mathcal{W}_{2}, \mathcal{X}_{p}\rangle \\
\bar{\gamma}_{-}=\frac{1}{m_{2}}\sum\limits_{q=1}^{m_{2}}\gamma_{-}^{q}=\frac{1}{m_{2}}\sum\limits_{q=1}^{m_{2}} y_{2}^{q}f_{1}(\mathcal{Y}_{q}) =\frac{1}{m_{2}}y_{2}^{q}\langle \mathcal{W}_{1}, \mathcal{Y}_{q}\rangle.
\end{cases}
\end{equation}
The margin variances are defined as follows:
\begin{equation}\label{3.2}
\begin{aligned}
\begin{cases}
\hat{\gamma}_{+}=\frac{1}{m_{1}}\sum\limits_{p=1}^{m_{1}}\left(\gamma_{+}^{p}-\bar{\gamma}_{+}\right)^{2}
=\frac{1}{m_{1}}\sum\limits_{p=1}^{m_{1}}\left(y_{1}^{p}\langle\mathcal{W}_{2}, \mathcal{X}_{p}\rangle-\bar{\gamma}_{+}\right)^{2}
=\left(\frac{m_{1}-1}{m_{1}^2}\right)\sum\limits_{p=1}^{m_{1}}\langle \mathcal{W}_{2},\mathcal{X}_{p}\rangle^{2}\\
%
%
\hat{\gamma}_{-}=\frac{1}{m_{2}}\sum\limits_{q=1}^{m_{2}}\left(\gamma_{-}^{q}-\bar{\gamma}_{-}\right)^{2}
=\frac{1}{m_{2}}\sum\limits_{q=1}^{m_{2}}\left(y_{2}^{q}\langle\mathcal{W}_{1}, \mathcal{Y}_{q}\rangle-\bar{\gamma}_{-}\right)^{2}
=\left(\frac{m_{2}-1}{m_{2}^2}\right)\sum\limits_{q=1}^{m_{2}}\langle \mathcal{W}_{1}, \mathcal{Y}_{q}\rangle^{2}
\end{cases}
\end{aligned}
\end{equation}

To effectively capture complex data relationships and enhance classification performance, we assume the weights in subsequent classifiers form a tensor \(\mathcal{W}\). This tensor can be decomposed into a sum of $R$ rank-one tensors, as defined by the CP decomposition outlined in Definition \ref{2.4}, i.e.
\begin{equation}\label{3.3}
\mathcal{W} = \sum_{r=1}^{R} \mathbf{u}^{(1)}_r \otimes \mathbf{u}^{(2)}_r \otimes \cdots \otimes \mathbf{u}^{(M)}_r,
\end{equation}
where \(\mathbf {u}^{(j)}_r \in \mathbb{R}^{I_j}, j = 1, 2, \ldots, M\), and $M$ is the number of tensor modes. The weights tensor $\mathcal{W}$ generalizes the weight vector \(\w\) in SVMs, where \(\w\)
represents the normal vector perpendicular to the separating hyperplane. The tensor form extends traditional vector-based weights to capture higher-dimensional structural information.

For mode \(j\) ($j = 1, \ldots, M$) stack the rank-one components $\{\mathbf{u}^{(j)}_r\}_{r=1}^R$ into a matrix:
$$\mathbf{U}^{(j)} = [\mathbf{u}^{(j)}_1, \mathbf{u}^{(j)}_2, \ldots, \mathbf{u}^{(j)}_R]\in\mathbb{R}^{I_j\times R}, $$
The \( j \)-th matricization of \(\mathcal{W}\) (i.e., reshaping the tensor into a matrix by unfolding along the \( j \)-th mode) is given by:
\begin{equation}\label{3.4}
\mathbf{W}^{(j)} = \mathbf{U}^{(j)}\left(\mathbf{U}^{(M)} \odot \cdots \odot \mathbf{U}^{(j+1)} \odot \mathbf{U}^{(j-1)} \odot \cdots \odot \mathbf{U}^{(1)}\right)^{\top}
=\mathbf{U}^{(j)}\left(\mathbf{U}^{(-j)}\right)^{\top}.
\end{equation}

Similarly, the \( j \)-th matricization of samples \(\mathcal{X}\) is given by $$\mathbf{X}^{(j)} =\mathbf{V}^{(j)}\left(\mathbf{V}^{(-j)}\right)^{\top}.$$
For samples \(\mathcal{Y}\), the matricization follows the same structural form (consistent with $\mathcal{X}$). 

Further, consider the tensor inner product property:
\begin{equation}\label{3.5}
\langle \mathcal{W}, \mathcal{W} \rangle = \Tr\left(\mathbf{W}^{(j)} {\mathbf{W}^{(j)}}^{\top}\right) = \text{vec}\left(\mathbf{W}^{(j)}\right)^{\top} \text{vec}\left(\mathbf{W}^{(j)}\right),
\end{equation}
where \(\text{Tr}(\cdot)\) denotes the trace operator, and \(\text{vec}(\cdot)\) denotes vectorization.

Analogous to Equation (\ref{3.5}), for individual sample tensors $\mathcal{X}_{p}$ and $\mathcal{Y}_{q}$, their inner products with weight tensors $\mathcal{W}_{1},\mathcal{W}_{2}$ satisfy:
\begin{equation}\label{3.6}
\begin{aligned}
\langle \mathcal{W}_{1}, \mathcal{X}_{p}\rangle=\Tr\left({\mathbf W}_{1}^{(j)}{\mathbf X_{p}^{(j)}}^{\top}\right),\quad\langle\mathcal{W}_{1}, \mathcal{Y}_{q}\rangle=\Tr\left({\mathbf W}_{1}^{(j)}{{\mathbf Y}_{q}^{(j)}}^{\top}\right),\\
\langle \mathcal{W}_{2}, \mathcal{X}_{p}\rangle=\Tr\left({\mathbf W}_{2}^{(j)}{\mathbf X_{p}^{(j)}}^{\top}\right),\quad\langle\mathcal{W}_{2}, \mathcal{Y}_{q}\rangle=\Tr\left({\mathbf W}_{2}^{(j)}{{\mathbf Y}_{q}^{(j)}}^{\top}\right).
\end{aligned}
\end{equation}
where ${\mathbf W}_{1}^{(j)},{\mathbf W}_{2}^{(j)}$  are the \( j \)-th matricizations of $\mathcal{W}_{1},\mathcal{W}_{2}$, and $\mathbf X_{p}^{(j)},{\mathbf Y}_{q}^{(j)}$ are the \( j \)-th
matricizations of $\mathcal{X}_{p},\mathcal{Y}_{q}$.

To guarantee that the matrices in the LMD-NPSTM dual problem are nonsingular, we add a regularization term to maximize some margin. The formula is described as

\begin{equation}\label{3.7}
\frac{c_i}{2} \|\mathcal{W}_i\|_F^2, \quad i = 1, 2.
\end{equation}
By utilizing (\ref{3.7}), the LMD-NPSTM dual problem can be derived without any additional assumptions and modifications.

LMD-NPSTM seeks a couple of tensors $\mathcal{W}_{1}$ and $\mathcal{W}_{2}$ simultaneously maximizing the mean of positive and negative margin while minimizing the variance of margin, i.e., focusing on two optimization tasks that follow
\begin{equation}\label{3.8}
\begin{aligned}
\min_{\mathcal{W}_{1},\bm{\xi}_{2}}\quad&\frac{1}{2}\sum\limits_{p=1}^{m_1}\langle \mathcal{W}_{1}, \mathcal{X}_{p}\rangle^{2}+\frac{c_1}{2}\|\mathcal{W}_{1}\|_F^2+\lambda_1\hat{\gamma}_{-}-\lambda_3
\bar{\gamma}_{-}+c_3\mathbf{e}_2^\top \bm{\xi}_2\\
{\text{s.t.}}\quad&-\sum_{q=1}^{m_2}\langle\mathcal{W}_{1}, \mathcal{Y}_{q}\rangle+\xi_2^q\geq1, ~\xi_2^q\geq0,~q=1,\ldots, m_2
\end{aligned}
\end{equation}
and
\begin{equation}\label{3.9}
\begin{aligned}
\min_{\mathcal{W}_{2},\bm{\xi}_{1}}\quad&\frac{1}{2}\sum\limits_{q=1}^{m_2}\langle \mathcal{W}_{2}, \mathcal{Y}_{q}\rangle^{2}+\frac{c_2}{2}\|\mathcal{W}_{2}\|_F^2+\lambda_2\hat{\gamma}_{+}-\lambda_4\bar{\gamma}_{+}+c_4\bm e_1^\top\bm\xi_1\\
{\text{s.t.}}\quad&-\sum_{p=1}^{m_1}\langle\mathcal{W}_{2}, \mathcal{X}_{p}\rangle+\xi_1^p\geq1, ~\xi_1^p\geq0,~p=1,\ldots, m_1,
\end{aligned}
\end{equation}
where $\lambda_1,\lambda_2,\lambda_3,\lambda_4$ are hyperparameters, which balance margin variance, margin mean, and model complexity;  
$c_1$ and $c_2$ control regularization strength, while $c_3$ and $c_4$ weight the penalty on slack variables $\bm\xi_1$ and $\bm\xi_2$; $\mathbf{e}_1$ and $\mathbf{e}_2$ are all-ones vectors matching the dimensions of $\bm\xi_1$ and $\bm\xi_2$.

However, the optimization problems \eqref{3.8}-\eqref{3.9} are non-convex in the tensor parameters $\mathcal{W}_{1}$ and $\mathcal{W}_{2}$. To address this, we adopt an alternating optimization scheme. 
At each iteration, fix all tensor modes except the $j$-th, solving for ${\mathbf W}_{1}^{(j)},{\mathbf W}_{2}^{(j)}$, while keeping other modes constant.
More specifically, and taking (\ref{3.5}) and (\ref{3.6}) into consideration, at the iterations for the $j$-th mode we solve the following optimization problems:
\begin{equation}\label{3.10}
\begin{aligned}
\min_{{\mathbf W}_{1}^{(j)},\bm\xi_{2}}\quad&\frac{1}{2}\sum\limits_{p=1}^{m_1}\left(\Tr\left({\mathbf W}_{1}^{(j)}{\mathbf X_{p}^{(j)}}^{\top}\right) \right)^2+
\frac{c_1}{2}\Tr\left(\mathbf W_1^{(j)}{\mathbf W_1^{(j)}}^{\top}\right)+\lambda_1\hat{\gamma}_{-}-\lambda_3\bar{\gamma}_{-}+c_3\bm e_2^\top\bm \xi_2\\
{\text{s.t.}}\quad&-\sum_{q=1}^{m_2}\Tr\left({\mathbf W}_{1}^{(j)}{{\mathbf Y}_{q}^{(j)}}^{\top}\right)+\xi_2^q\geq1, ~\xi_2^q\geq0,~q=1,\ldots, m_2
\end{aligned}
\end{equation}
and
\begin{equation}\label{3.11}
\begin{aligned}
\min_{{\mathbf W}_{2}^{(j)},\bm\xi_{1}}\quad&\frac{1}{2}\sum\limits_{q=1}^{m_2}\left(\Tr\left({\mathbf W}_{2}^{(j)}{\mathbf Y_{q}^{(j)}}^{\top}\right)\right)^2+
\frac{c_2}{2}\Tr\left(\mathbf W_2^{(j)}{\mathbf W_2^{(j)}}^{\top}\right)+\lambda_2\hat{\gamma}_{-}-\lambda_4\bar{\gamma}_{-}+c_4\bm e_1^\top\bm \xi_1\\
{\text{s.t.}}\quad&-\sum_{p=1}^{m_1}\Tr\left({\mathbf W}_{2}^{(j)}{{\mathbf X}_{p}^{(j)}}^{\top}\right)+\xi_1^p\geq1, ~\xi_1^p\geq0,~p=1,\ldots, m_1.
\end{aligned}
\end{equation}

Under our assumption/constraint that the tensor $\mathcal{W}_{1}$ and $\mathcal{W}_{2}$ are written as a sum of rank one tensors as in (\ref{3.3}), we replace ${\mathbf W}_{1}^{(j)}$ and ${\mathbf W}_{2}^{(j)}$ in the above equation by taking (\ref{3.7}) into consideration. We should mention here that the initial values for the rank-one tensors $\bm u_1$ and $\bm u_2$, and subsequently for their matricized forms $\U_{1}^{(j)}$ and $\U_{2}^{(j)}$, are randomly chosen. Then (\ref{3.10}) and (\ref{3.11}) are rewritten as
\begin{equation}\label{3.12}
\begin{aligned}
\min_{{\U}_{1}^{(j)},\bm \xi_{2}}\quad&\frac{1}{2}\sum\limits_{p=1}^{m_1}
\left(\Tr\left({\U}_{1}^{(j)}\left({\U}_{1}^{(-j)}\right)^{\top}{\X_{p}^{(j)}}^{\top}\right)\right)^2+
\frac{c_1}{2}\Tr\left({\U}_{1}^{(j)}\left({\U}_{1}^{(-j)}\right)^{\top}\left({\U}_{1}^{(-j)}\right)\left({\U}_{1}^{(j)}\right)^{\top}\right)\\
&+\lambda_1\hat{\gamma}_{-}-\lambda_3\bar{\gamma}_{-}
+c_3\bm e_2^\top\bm \xi_2\\
{\rm s.t.}\quad&-\sum_{q=1}^{m_2}\Tr\left({\U}_{1}^{(j)}\left({\U}_{1}^{(-j)}\right)^{\top}{{\Y}_{q}^{(j)}}^{\top}\right)+\xi_2^q\geq1, ~\xi_2^q\geq0,~q=1,\ldots, m_2.
\end{aligned}
\end{equation}


and
\begin{equation}\label{3.13}
\begin{aligned}
\min_{{\U}_{2}^{(j)},\bm\xi_{1}}\quad&\frac{1}{2}\sum\limits_{q=1}^{m_2}
\left(\Tr\left({\U}_{2}^{(j)}\left({\U}_{2}^{(-j)}\right)^{\top}{\Y_{q}^{(j)}}^{\top}\right)\right)^2+
\frac{c_2}{2}\Tr\left({\U}_{2}^{(j)}\left({\U}_{2}^{(-j)}\right)^{\top}\left({\U}_{2}^{(-j)}\right)\left({\U}_{2}^{(j)}\right)^{\top}\right)\\
&+\lambda_2\hat{\gamma}_{-}-\lambda_4\bar{\gamma}_{-}
+c_4\bm e_1^\top\bm \xi_1\\
{\rm s.t.}\quad&-\sum_{p=1}^{m_1}\Tr\left({\U}_{2}^{(j)}\left({\U}_{2}^{(-j)}\right)^{\top}{{\X}_{p}^{(j)}}^{\top}\right)+\xi_1^p\geq1, ~\xi_1^p\geq0,~p=1,\ldots, m_1.
\end{aligned}
\end{equation}

To simplify the trace operations in \eqref{3.12}-\eqref{3.13}, we introduce a change of variables leveraging the positive definiteness of factor matrices. Define $\A={{\U}_{1}^{(-j)}}^{\top}{\U}_{1}^{(-j)}$, which is a positive matrix. Introduce the transformed factor matrix:
$\tilde{{\U}}_{1}^{(j)}={\U}_{1}^{(j)}\A^{\frac{1}{2}}$. The regularization trace in \eqref{3.12} becomes:
\begin{equation}\label{3.14}
\Tr\left({\U}_{1}^{(j)}\left({\U}_{1}^{(-j)}\right)^{\top}\left({\U}_{1}^{(-j)}\right)
\left({\U}_{1}^{(j)}\right)^{\top}\right)
=\Tr\left(\tilde{{\U}}_{1}^{(j)}\left(\tilde{{\U}}_{1}^{(j)}\right)^{\top}\right)=
{\rm vec}\left(\tilde{{\U}}_{1}^{(j)}\right)^{\top}{\rm vec}\left(\tilde{{\U}}_{1}^{(j)}\right)
\end{equation}
Define transformed data matrices for mode $j$: $\tilde{\X}_p^{(j)}=\X_{p}^{(j)}{\U}_{1}^{(-j)}\A^{-\frac{1}{2}},~
\tilde{\Y}_q^{(j)}=\Y_{q}^{(j)}{\U}_{1}^{(-j)}\A^{-\frac{1}{2}}$, The data fidelity trace in \eqref{3.12} simplifies to:
\begin{equation}\label{3.15}
\Tr\left({\U}_{1}^{(j)}\left({\U}_{1}^{(-j)}\right)^{\top}{\X_{p}^{(j)}}^{\top}\right)
=\Tr\left(\tilde{{\U}}_{1}^{(j)}\left(\tilde{\X}_p^{(j)}\right)^{\top}\right)
={\rm vec}\left(\tilde{{\U}}_{1}^{(j)}\right)^{\top}{\rm vec}\left(\tilde{\X}_p^{(j)}\right)
\end{equation}
\begin{equation}\label{3.16}
\Tr\left({\U}_{1}^{(j)}\left({\U}_{1}^{(-j)}\right)^{\top}{\Y_{q}^{(j)}}^{\top}\right)
=\Tr\left(\tilde{{\U}}_{1}^{(j)}\left(\tilde{\Y}_q^{(j)}\right)^{\top}\right)=
{\rm vec}\left(\tilde{{\U}}_{1}^{(j)}\right)^{\top}{\rm vec}\left(\tilde{\Y}_q^{(j)}\right)
\end{equation}
Then, the margin mean/variance (\ref{3.1})-(\ref{3.2}), with the transformed variables \eqref{3.14}-\eqref{3.16}, the subproblem (\ref{3.12}) becomes 
\begin{equation} \label{3.17}
\begin{aligned}
\min_{{\text{vec}}\left(\tilde{{\U}}_{1}^{(j)}\right),\bm \xi_{2}}
\quad&\frac{1}{2}\sum\limits_{p=1}^{m_1}\left({\rm vec}
\left(\tilde{{\U}}_{1}^{(j)}\right)^{\top}{\rm vec}\left(\tilde{\X}_p^{(j)}\right)\right)^2-\frac{\lambda_3}{m_2}\sum\limits_{q=1}^{m_2}\y_2^{q}{\rm vec}\left(\tilde{{\U}}_{1}^{(j)}\right)^{\top}{\rm vec}\left(\tilde{\Y}_q^{(j)}\right)+c_3\bm e_2^\top\bm \xi_2\\
&+\frac{c_1}{2}\text{vec}\left(\tilde{{\U}}_{1}^{(j)}\right)^{\top}{\rm vec}\left(\tilde{{\U}}_{1}^{(j)}\right)
+\left(\frac{\lambda_1(m_{2}-1)}{m_{2}^2}\right)\sum\limits_{q=1}^{m_2}\left({\rm vec}\left(\tilde{{\U}}_{1}^{(j)}\right)^{\top}{\rm vec}\left(\tilde{\Y}_q^{(j)}\right)\right)^2\\
{\rm s.t.}\quad&-{\rm vec}\left(\tilde{{\U}}_{1}^{(j)}\right)^{\top}{\rm vec}\left(\tilde{\Y}_q^{(j)}\right)+\xi_2^q\geq1, ~\xi_2^q \geq0,~q=1,\ldots, m_2.
\end{aligned}
\end{equation}

Similarly, Let us define $\B=\left({\U}_{2}^{(-j)}\right)^{\top}{\U}_{2}^{(-j)},~ \tilde{{\U}}_{2}^{(j)}={\U}_{2}^{(j)}\B^{\frac{1}{2}}$, where $\B$ is a positive matrix. Then,
\begin{equation}\label{3.18}
\Tr\left({\U}_{2}^{(j)}\left({\U}_{2}^{(-j)}\right)^{\top}\left({\U}_{2}^{(-j)}\right)\left({\U}_{2}^{(j)}\right)^{\top}\right)
=\Tr\left(\tilde{{\U}}_{2}^{(j)}\left(\tilde{{\U}}_{2}^{(j)}\right)^{\top}\right)=
{\rm vec}\left(\tilde{{\U}}_{2}^{(j)}\right)^{\top}{\rm vec}\left(\tilde{{\U}}_{2}^{(j)}\right)
\end{equation}
The transformed data matrices are: $\bar{\Y}_q^{(j)}=\Y_{q}^{(j)}{\U}_{2}^{(-j)}\B^{-\frac{1}{2}},
~\bar{\X}_p^{(j)}=\X_{p}^{(j)}{\U}_{2}^{(-j)}\B^{-\frac{1}{2}}$,
we have
\begin{equation}\label{3.19}
\Tr\left({\U}_{2}^{(j)}\left({\U}_{2}^{(-j)}\right)^{\top}{\Y_{q}^{(j)}}^{\top}\right)
=\Tr\left(\tilde{{\U}}_{2}^{(j)}\left(\bar{\Y}_q^{(j)}\right)^{\top}\right)={\rm vec}\left(\tilde{{\U}}_{2}^{(j)}\right)^{\top}{\rm vec}\left(\bar{\Y}_q^{(j)}\right)
\end{equation}
\begin{equation}\label{3.20}
\Tr\left({\U}_{2}^{(j)}\left({\U}_{2}^{(-j)}\right)^{\top}{\X_{p}^{(j)}}^{\top}\right)
=\Tr\left(\tilde{{\U}}_{2}^{(j)}\left(\tilde{\X}_p^{(j)}\right)^{\top}\right)=
{\rm vec}\left(\tilde{{\U}}_{2}^{(j)}\right)^{\top}{\rm vec}\left(\tilde{\X}_p^{(j)}\right)
\end{equation}
Then, combining (\ref{3.1}) with (\ref{3.2}), (\ref{3.13}) is written as
\begin{equation}\label{3.21}
\begin{aligned}
\min_{{\rm vec}\left(\tilde{{\U}}_{2}^{(j)}\right),\bm\xi_{1}}
\quad&\frac{1}{2}\sum\limits_{q=1}^{m_2}\left({\rm vec}\left(\tilde{{\U}}_{2}^{(j)}\right)^{\top}
{\rm vec}\left(\tilde{\Y}_q^{(j)}\right)\right)^2
-\frac{\lambda_4}{m_1}\sum\limits_{p=1}^{m_1}\y_1^{p}
{\rm vec}\left(\tilde{{\U}}_{2}^{(j)}\right)^{\top}{\rm vec}\left(\tilde{\X}_p^{(j)}\right)+c_4 \bm e_1^\top\bm \xi_1\\
&+\frac{c_2}{2}\left(\tilde{{\U}}_{2}^{(j)}\right)^{\top}{\rm vec}\left(\tilde{{\U}}_{2}^{(j)}\right)+\left(\frac{\lambda_2(m_{1}-1)}{m_{1}^2}\right)
\sum\limits_{p=1}^{m_1}\left({\rm vec}(\tilde{{\U}}_{2}^{(j)})^{\top}{\rm vec}\left(\tilde{\X}_p^{(j)}\right)\right)^2\\
{\rm s.t.}\quad&-{\rm vec}\left(\tilde{{\U}}_{2}^{(j)}\right)^{\top}{\rm vec}\left(\tilde{\X}_p^{(j)}\right)+\xi_1^p\geq1, ~\xi_1^p\geq0,~p=1,\ldots, m_1.
\end{aligned}
\end{equation}

\begin{theorem}\label{th4.1}
The optimal solution ${{\rm vec}\left(\tilde{{\U}}_{1}^{(j)}\right)}^\ast$ and ${\rm vec}\left(\tilde{{\U}}_{2}^{(j)}\right)^\ast$ of (\ref{3.17}) and (\ref{3.21})  can be expressed succinctly as ${\rm vec}\left(\tilde{{\U}}_{1}^{(j)}\right)=\bm V^{(j)}\bm\beta_1$ and ${\rm vec}\left(\tilde{{\U}}_{2}^{(j)}\right)=\bm V^{(j)}\bm\beta_2$, where $\beta_1,\beta_2\in \mathbb{R}^m$ are coeffcient vectors.
\end{theorem}

\proof To analyze the objective function, we decompose the vectorized factor matrix ${\rm vec}\left(\tilde{{\U}}_{1}^{(j)}\right)$ using the span of sample modes. By the projection theorem,
${\rm vec}\left(\tilde{{\U}}_{1}^{(j)}\right)$ can be decomposed into a part that lives in the span of ${\rm vec}\left(\bm V_l^{(j)}\right)$ and an orthogonal part vector, i.e.,
\begin{equation*}
{\rm vec}\left(\tilde{{\U}}_{1}^{(j)}\right)=\sum\limits_{l=1}^m\bm\beta_1^l{\rm vec}\left(\bm V_l^{(j)}\right)+\bm\eta_1=\bm V^{(j)}\bm\beta_1+\bm\eta_1,~ l=1,2,\ldots,m,
\end{equation*}
where $\bm\eta_1$ is a vector and $m=m_1+m_2$, which satisfies $\left(\bm V^{(j)}\right)^{\top}\bm\eta_1=\0$ \cite{zz14}.

The objective function of (\ref{3.17}) includes a quadratic term in ${\rm vec}\left(\tilde{{\U}}_{1}^{(j)}\right)$:
\begin{equation}\label{3.22}
\begin{aligned}
{\rm vec}\left(\tilde{{\U}}_{1}^{(j)}\right)^{\top}{\rm vec}\left(\tilde{{\U}}_{1}^{(j)}\right)
&=\left(\bm V^{(j)}\bm\beta_1+\bm\eta_1\right)^{\top}\left(\bm V^{(j)}\bm\beta_1+\bm\eta_1\right)\\
&=\bm\beta_1^{\top}\left(\bm V^{(j)}\right)^{\top}\bm V^{(j)}\bm\beta_1+\bm\eta_1^{\top}\bm\eta_1\\
\end{aligned}
\end{equation}
By the orthogonality $\left(\bm V^{(j)}\right)^{\top}\bm\eta_1=\0$, the cross term vanishes. Define $\K=\left(\bm V^{(j)}\right)^{\top}\bm V^{(j)}$, then $${\rm vec}\left(\tilde{{\U}}_{1}^{(j)}\right)^{\top}{\rm vec}\left(\tilde{{\U}}_{1}^{(j)}\right)=\bm\beta_1^{\top}\K\bm\beta_1+\|\bm\eta_1\|^2.$$
Since $\|\bm\eta_1\|^2\geq0$, the quadratic term is minimized when ${\bm\eta}_1=\0$. This implies
${\rm vec}\left(\tilde{{\U}}_{1}^{(j)}\right)=\bm V^{(j)}\bm\beta_1$.

A parallel decomposition applies to ${\rm vec}\left(\tilde{{\U}}_{2}^{(j)}\right)$, 
${\rm vec}\left(\tilde{{\U}}_{2}^{(j)}\right)=\bm V^{(j)}\bm\beta_2+\bm\eta_2$, with $\left(\bm V^{(j)}\right)^{\top}\bm\eta_2=\0$, leading to the same conclusion: ${\bm\eta}_2=\0$ minimizes the quadratic term. 

So, setting $\bm\eta_i~ (i=1,2)=\0$ does not affect the other terms, but strictly reduces the fourth term of the objective function. Therefore,
Theorem \ref{th4.1} is accepted.

Based on Theorem \ref{th4.1}, the transformation of the original problem into the Wolfe dual problem can be derived.

\begin{theorem}\label{th4.2}
The original problem (\ref{3.8}) can be transformed into the Wolfe dual problem as follows:
\begin{equation*}
\begin{aligned}
\max_{\bm{\alpha}_1}\quad&-\frac{1}{2}\bm{\alpha}_1^{\top} \H_1\bm{\alpha}_1+\left(\frac{\lambda_3}{m_2}\H_1\y_2+\e_2\right)^{\top}\bm{\alpha}_1\\
{\rm s.t.}\quad&\0\leq\bm{\alpha}_1 \leq c_3\bm {s}_3
\end{aligned}
\end{equation*}
\noindent where $\bm {H}_1=\bm K_2\bm {G}_1^{-1}  \bm {K}_2, \G_1=\K_1^{\top}\K_1+c_1\K+2\lambda_1\left(\frac{m_2-1}{m_2^2}\right)\K_2^{\top}\K_2$, $\bm \alpha_1$ is a Lagrangian multipliers vector and $\bm\beta_1$ has the following representation:
\begin{equation*}
\bm \beta_1= \G_1^{-1} \left(\frac{\lambda_3}{m_2}\K_2^{\top}y_2-\K_2^{\top}\bm\alpha_1\right).
\end{equation*}
\end{theorem}
\proof According to \eqref{3.22} from Theorem \ref{th4.1}, we have
\begin{equation}\label{3.23}
\begin{aligned}
\sum\limits_{p=1}^{m_1}\left({\rm vec}\left(\tilde{{\U}}_{1}^{(j)}\right)^{\top}{\rm vec}\left(\tilde{\X}_p^{(j)}\right)\right)^2
&=\left(\left(\tilde{\X}_p^{(j)}\right)^{\top}{\rm vec}\left(\tilde{{\U}}_{1}^{(j)}\right)\right)^{\top}\left(\left(\tilde{\X}_p^{(j)}\right)^{\top}{\rm vec}\left(\tilde{{\U}}_{1}^{(j)}\right)\right)\\
&=\left(\left(\tilde{\X}_p^{(j)}\right)^{\top}\bm V^{(j)}\bm\beta_1\right)^{\top}\left(\left(\tilde{\X}_p^{(j)}\right)^{\top}\bm V^{(j)}\bm\beta_1\right)\\
&=\bm\beta_1^{\top}\left(\bm V^{(j)}\right)^{\top}\tilde{\X}_p^{(j)}\left(\tilde{\X}_p^{(j)}\right)^{\top}\bm V^{(j)}\bm\beta_1\\
&=\bm\beta_1^{\top}\K_1^{\top}\K_1\bm\beta_1,
\end{aligned}
\end{equation}
where $\K_1=\left(\tilde{\X}_p^{(j)}\right)^{\top}\bm V^{(j)}$,
\begin{equation}\label{3.24}
\begin{aligned}
\sum\limits_{q=1}^{m_2}\y_2^{q}{\rm vec}\left(\tilde{{\U}}_{1}^{(j)}\right)^{\top}{\rm vec}\left(\tilde{\Y}_q^{(j)}\right)
&=\y_2^{\top}\left(\tilde{\Y}_q^{(j)}\right)^{\top}{\rm vec}\left(\tilde{{\U}}_{1}^{(j)}\right)\\
&=\y_2^{\top}\left(\tilde{\Y}_q^{(j)}\right)^{\top}\bm V^{(j)}\bm\beta_1\\
&=\y_2^{\top}\bm M_1\bm\beta_1,
\end{aligned}
\end{equation}
where $\M_1=\left(\tilde{\Y}_q^{(j)}\right)^{\top}\bm V^{(j)}$,
and
\begin{equation}\label{3.25}
\begin{aligned}
\sum\limits_{q=1}^{m_2}\left({\rm vec}\left(\tilde{{\U}}_{1}^{(j)}\right)^{\top}{\rm vec}\left(\tilde{\Y}_q^{(j)}\right)\right)^2
&=\left(\left(\tilde{\Y}_q^{(j)}\right)^{\top}{\rm vec}\left(\tilde{{\U}}_{1}^{(j)}\right)\right)^{\top}\left(\left(\tilde{\Y}_q^{(j)}\right)^{\top}{\rm vec}\left(\tilde{{\U}}_{1}^{(j)}\right)\right)\\
&=\left(\left(\tilde{\Y}_q^{(j)}\right)^{\top}\bm V^{(j)}\bm\beta_1\right)^{\top}\left(\left(\tilde{\Y}_q^{(j)}\right)^{\top}\bm V^{(j)}\bm\beta_1\right)\\
&=\bm\beta_1^{\top}\left(\bm V^{(j)}\right)^{\top}\tilde{\Y}_q^{(j)}\left(\tilde{\Y}_q^{(j)}\right)^{\top}\bm V^{(j)}\bm\beta_1\\
&=\bm\beta_1\M_1^{\top}\M_1\bm\beta_1
\end{aligned}
\end{equation}
Substituting (\ref{3.23}), (\ref{3.24}) and (\ref{3.25}) into (\ref{3.17}), we can obtain the following final matrix form
\begin{equation}\label{3.26}
\begin{aligned}
\min_{\bm\beta_1,\bm\xi_{2}}\quad&\frac{1}{2}\bm\beta_1^{\top}\K_1^{\top}\K_1\bm\beta_1+\frac{c_1}{2}\bm\beta_1^{\top}\K\bm\beta_1
+\left(\frac{\lambda_1(m_{2}-1)}{m_{2}^2}\right)\bm\beta_1^{\top}\M_1^{\top}\M_1\bm\beta_1-\frac{\lambda_3}{m_2}\y_2^{\top}\M_1\bm\beta_1+c_3
\bm e_2^\top\bm\xi_2\\
{\rm s.t.}\quad&-\M_1\bm\beta_1+\bm\xi_2\geq \e_2, ~\bm\xi_2\geq\0.
\end{aligned}
\end{equation}
And (\ref{3.26}) can be written in a more concise form as shown below:
\begin{equation}\label{3.27}
\begin{aligned}
\min_{\bm\beta_1,\bm\xi_{2}}\quad&\frac{1}{2}\bm\beta_1^{\top}\G_1\bm\beta_1-\frac{\lambda_3}{m_2}\y_2^{\top}\M_1\bm\beta_1+c_3
\bm e_2^\top\bm\xi_2\\
{\rm s.t.}\quad&-\M_1\bm\beta_1+\bm\xi_2\geq \e_2, ~\bm\xi_2\geq\0,
\end{aligned}
\end{equation}
where $\G_1=\K_1^{\top}\K_1+c_1\K+\left(\frac{2\lambda_1(m_{2}-1)}{m_{2}^2}\right)\M_1^{\top}\M_1$, and is a symmetric nonnegative definite matrix.

The Lagrangian function of the optimization problem in (\ref{3.27}) is presented as
\begin{equation}\label{3.28}
L_1(\bm\beta_1,\bm\xi_2,\bm\alpha_1,\bm\delta_1)=\frac{1}{2}\bm\beta_1^{\top}\G_1\bm\beta_1-\frac{\lambda_3}{m_2}\y_2^{\top}\M_1\bm\beta_1+c_3
\bm e_2^\top\bm\xi_2-\bm\alpha_1^{\top}(-\M_1\bm\beta_1+\bm\xi_2-\e_2)-\bm\delta_1^{\top}\bm\xi_2,
\end{equation}
where $\bm\alpha_1,\bm\delta_1\in\mathbb{R}^{m_2}$ are Lagrangian multipliers vectors.

According to the dual theorem and Karush-Kuhn-Tucker (K.K.T.), the minimum of the Lagrangian function in (\ref{3.28}) with respect to  $\bm\beta_1, \bm\xi_2 $ equals the maximum of the function with respect to $\bm\alpha_1$. By satisfying the necessary
conditions for the optimal solution for the Lagrange function, i.e.,
$\frac{\partial L_1}{\partial \bm\beta_1}=\frac{\partial L_1}{\partial \bm\xi_2}=\0$, we obtain the following formulas:
\begin{equation}\label{3.29}
\begin{aligned}
\G_1\bm\beta_1=\frac{\lambda_3}{m_2}\y_2^{\top}\M_1-\bm\alpha_1^{\top}\M_1\\
c_3\bm e_2-\bm\alpha_1-\bm\delta_1=\0\Rightarrow\0\leq\bm\alpha_1\leq c_3\bm e_2
\end{aligned}
\end{equation}
Since $\G_1$ is nonsingular, $\bm\beta_1$ can be deduced from (\ref{3.29}) as mentioned below:
\begin{equation}\label{3.30}
\bm\beta_1=\G_1^{-1}\left(\frac{\lambda_3}{m_2}\y_2^{\top}\M_1-\bm\alpha_1^{\top}\M_1\right).
\end{equation}
Submitting (\ref{3.29}) into the Lagrangian function (\ref{3.28}), the Wolfe dual form of the model (\ref{3.27}) is presented as follows:
\begin{equation}\label{3.31}
\begin{aligned}
\max_{\bm\alpha_1}\quad&-\frac{1}{2}\bm\alpha_1^{\top}\bm H_1\bm\alpha_1+\left(\frac{\lambda_3}{m_2}\bm H_1\y_2+\bm e_2\right)^{\top}\bm\alpha_1\\
{\rm s.t.}\quad&\0\leq\bm\alpha_1\leq c_3\bm e_2.
\end{aligned}
\end{equation}
where $\bm H_1=\M_1\bm G_1^{-1}\M_1^{\top}$.

Similarly,
we derive the the Wolfe dual form of model (\ref{3.21}) using Theorem \ref{th4.1} as indicated bellow:
\begin{equation}\label{3.32}
\begin{aligned}
\max_{\bm\alpha_2}\quad&-\frac{1}{2}\bm\alpha_2^{\top}\bm H_2\bm\alpha_2+\left(\frac{\lambda_4}{m_1}\bm H_2\y_1\bm e_1\right)^{\top}\bm\alpha_2\\
{\rm s.t.}\quad&\0\leq\bm\alpha_2\leq c_4\bm e_1.
\end{aligned}
\end{equation}
where $\bm\alpha_2$ is a nonnegative Lagrangian multipliers vectors, and $\bm \H_2=\bm M_2\bm G_2^{-1}\bm M_2^{\top}$ and $\bm\beta_2$ is expressed as
\begin{equation*}
\bm\beta_2=\G_2^{-1}\left(\frac{\lambda_4}{m_1}\y_1^{\top}\bm M_2-\bm\alpha_2^{\top}\bm M_2\right)
\end{equation*}
where
$\G_2=\bm K_2^{\top}\bm K_2+c_2\K+\left(\frac{2\lambda_2(m_{1}-1)}{m_{1}^2}\right)\bm M_2^{\top}\bm M_2$,  and is a symmetric nonnegative definite matrix, and $\bm K_2=\left(\bar{\Y}_q^{(j)}\right)^{\top}\bm V^{(j)}, \bm M_2=\left(\bar{\X}_p^{(j)}\right)^{\top}\bm V^{(j)}$.

\begin{theorem}
By utilizing Theorem \ref{th4.2}, we obtain the Wolfe dual problems and the parameter vectors \( \bm\beta_1 \) and \( \bm\beta_2 \) related to the optimal hyperplanes. Then, the expression for the decision function of LMD-NPSTM can be constructed as follows:
$$
f(\mathcal{X}) = \arg \min_{i=1,2} \frac{|\langle\mathcal{W}_i,\mathcal{X}\rangle|}{\langle\mathcal{W}_i,\mathcal{W}_i\rangle}
=\frac{\left|{\rm vec}\left(\tilde{{\U}}_{i}^{(j)}\right)^{\top}{\rm vec}\left(\tilde{{\X}}^{(j)}\right)\right|}{{\rm vec}\left(\tilde{{\U}}_{i}^{(j)}\right)^{\top}{\rm vec}\left(\tilde{{\U}}_{i}^{(j)}\right)}
=\frac{|\K_i\bm\beta_i|}{\sqrt{\bm\beta_i^{\top} \K \bm\beta_i}}
$$
and
$$
f(\mathcal{Y}) = \arg \min_{i=1,2} \frac{|\langle\mathcal{W}_i,\mathcal{Y}\rangle|}{\langle\mathcal{W}_i,\mathcal{W}_i\rangle}
=\frac{\left|{\rm vec}\left(\tilde{{\U}}_{i}^{(j)}\right)^{\top}{\rm vec}\left(\tilde{{\Y}}^{(j)}\right)\right|}{{\rm vec}\left(\tilde{{\U}}_{i}^{(j)}\right)^{\top}{\rm vec}\left(\tilde{{\U}}_{i}^{(j)}\right)}
=\frac{|\M_i\bm\beta_i|}{\sqrt{\bm\beta_i^{\top} \K \bm\beta_i}}.
$$
\end{theorem}
\proof Obviously, the decision function is defined as the class of the hyperplane that is closer to the input point, and the concrete proof is similar to the traditional TSVMs \cite{krs07}.

\setcounter{proposition}{0} \setcounter{theorem}{0}
\setcounter{lemma}{0} \setcounter{corollary}{0}
\setcounter{equation}{0}
\setcounter{definition}{0}
\subsection{Algorithm and pseudocode}
In this section, we address problems (\ref{3.3}) and (\ref{3.9}) by proposing an alternating projection algorithm. For subproblems (\ref{3.10}) and (\ref{3.11}), we leverage its dual counterpart (\ref{3.31}) and (\ref{3.32}) to derive an efficient solution method.  The subsequent content is structured as follows:  we first delineate the iterative procedure of the algorithm (summarized concisely in Algorithm \ref{alg1} and then rigorously establish its convergence properties through theoretical analysis.

\begin{algorithm}[ht]
\caption {Alternating projection for LMD-NPSTM}
\label{alg1}
\algnewcommand\algorithmicinput{\textbf{Input:}}
\algnewcommand\algorithmicoutput{\textbf{Output:}}
\algnewcommand\Input{\item[\algorithmicinput]}
\algnewcommand\Output{\item[\algorithmicoutput]}
\setstretch{1.1}
\begin{algorithmic}[1]
\Input The set of training tensors \( \mathcal{X}_p|_{p=1}^{m_1},~\mathcal{Y}_q|_{q=1}^{m_2}\in \mathbb{R}^{I_1 \times I_2 \times \cdots \times I_M} \) and their corresponding labels \( y_i \in \{ +1, -1 \} \), N=5000.
\Output The parameters of the classification tensors \( \mathcal{W}_1 \) and \( \mathcal{W}_2 \). \STATE Initialize \( \mathcal{W}_1 \) and  \( \mathcal{W}_2 \) written as a sum of random rank-one tensors and $k=0$.
\WHILE {\( \|\mathcal{W}_i^{(k)} - \mathcal{W}_i^{(k-1)}\| / \|\mathcal{W}_i^{(k-1)}\| > \epsilon \) $(i=1,2)$ or $k<N$}
  \FOR {$ j = 1, \ldots, M $} (number of modes)
    \STATE Update variables $\bm\alpha_1$ in (\ref{3.31}) and $\bm\alpha_2$ in (\ref{3.32}) using the QPP toolkit.
    \STATE $\bm\beta_1\leftarrow\bm\alpha_1$, $\bm\beta_2\leftarrow\bm\alpha_2$.
    \STATE ${\rm vec}\left(\tilde{{\U}}_{1}^{(j)}\right)\leftarrow\bm\beta_1, {\rm vec}\left(\tilde{{\U}}_{2}^{(j)}\right)\leftarrow\bm\beta_2$.
   \STATE Calculate ${\U}_{1}^{(j)}$ and ${\U}_{1}^{(-j)}$ by (\ref{3.14}), ${\U}_{2}^{(j)}$ and ${\U}_{2}^{(-j)}$ by (\ref{3.18}).
\ENDFOR \\
$\mathbf{W}_1^{(j)}(k) \leftarrow \mathbf{U}_1^{(j)}\left(\mathbf{U}_1^{(-j)}\right)^{\top},\mathbf{W}_2^{(j)}(k) \leftarrow \mathbf{U}_2^{(j)}\left(\mathbf{U}_2^{(-j)}\right)^{\top} $.\\
Finish calculation of $\mathcal{W}_i^{(k)}$  ($i=1,2$).
\ENDWHILE
\end{algorithmic}
\end{algorithm}

The convergence of the algorithm is given below.

\begin{theorem}
Assume that $\left\{\left(\U_i^{(j)}, \bm\xi_{-i}^{(j)}\right)\right\}_{j=1}^M$ are the sequences generated in Algorithm \ref{alg1}. Here, $\bm\xi_{-i}^{(j)}$ represents the solution to the $j$-th subproblem of problem (\ref{3.8}) and (\ref{3.9}), respectively.  Then, the sequences $\left\{ f_i\left(\U_i^{(j)}, \bm\xi_{-i}^{(j)}\right)\right\}_{j=1}^M$ are each monotonically non-increasing and converge to limit points.
\end{theorem}
\proof We define the index set \( I = \{1, 2\} \), where for each \( i \in I \), the notation \( -i \)  denotes the complementary element in \( I \) (i.e., \( -i = 2 \) if \( i = 1 \), and vice versa).  We aim to minimize problems (\ref{3.8}) and (\ref{3.9}). Specifically, we consider minimizing certain functions $f_i({\U_i^{(1)},\ldots,\U_i^{(M)}},\bm\xi_{-i})$ that is of the form $f_i:\mathrm{R}^{I_1\times R}\times\cdots\times\mathrm{R}^{I_M\times R}\times\mathrm{R}^{m_i}\rightarrow\mathrm{R}$, since ${\U}_{i}^{(j)}\in\mathrm{R}^{I_j\times R}$.
The alternating optimization procedure can be written as
\begin{equation*}
g_i^j\!\left({\U_i^{(j)}}^\ast,{\bm\xi_{-i}^{(j)}}^\ast\right) \triangleq \min_{\U_i^{(j)}\in \mathbb{R}^{I_j\times R},\bm\xi_{-i}^{(j)}} f_i\!\left(\U_i^{(j)}, \overline{\{ \U_i^{(j)} \}}, \bm \xi_{-i}\right),
\end{equation*}
where $\overline{\left\{ \U_i^{(j)}\right \}} = \left\{ \U_i^{(1)}, \ldots, \U_i^{(j-1)}, \U_i^{(j+1)}, \ldots, \U_i^{(M)}\right\}$ denotes the set of all variables except $\U_i^{(j)}$. The superscript $\ast$ indicates the optimal solutions at iteration $j$.

In each iteration $k$, the function $g_i^j$ is computed using $\U_i^{(1)}, \ldots, \U_i^{(j-1)}$ from the current iteration $k$ and $\U_i^{(j+1)}, \ldots, \U_i^{(M)}$ obtained the previous iteration $k-1$.

Given an initialization of $\left\{ \U_i^{(j)} \right\}_{j=1}^M$, the alternating projection generates a sequence of $\left\{\left( \U_i^{(j)}(k), \bm\xi_{-i}^{(j)}(k)\right)\right\}_{j=1}^M$ such that
\begin{equation*}
g_i^j\!\left({\U_i^{(j)}}^\ast(k), {\bm\xi_{-i}^{(j)}}^\ast(k)\right) = \min_{\U_i^{(j)} \in \mathbb{R}^{I_j \times m},\bm\xi_{-i}} f_i\left(\U_i^{(j)}, \left\{ \U_i^{(l)}(k) \right\}_{l=1}^{j-1}, \left\{ \U_i^{(l)}(k-1) \right\}_{l=j+1}^{M},\bm\xi_{-i}\right).
\end{equation*}

By construction, each subproblem $g_i^j$ is a minimization step that guarantees:
$$f_i\!\left({\U_i^{(j)}}^\ast(k), \overline{\{ \U_i^{(j)} \}}, {\bm \xi_{-i}^{(j)}}^\ast(k)\right)\leq f_i\!\left({\U_i^{(j)}}^\ast(k-1), \overline{\{ \U_i^{(j)} \}}, {\bm \xi_{-i}^{(j)}}^\ast(k-1)\right),$$
Specifically, for any $j$ and $k$:
$$g_i^j\!\left({\U_i^{(j)}}^\ast(k), {\bm \xi_{-i}^{(j)}}^\ast(k)\right)\leq g_i^{j-1}\left({\U_i^{(j-1)}}^\ast(k), {\bm \xi_{-i}^{(j)}}^\ast(k)\right).$$
Therefore, the following holds:
\begin{equation*}
\begin{aligned}
\bar{\gamma_i} &= g_i^1\left({\U_i^{(1)}}^\ast(1),{\bm\xi_{-i}^{(1)}}^\ast(1)\right)\geq g_i^2\left({\U_i^{(2)}}^\ast(1), {\bm\xi_{-i}^{(1)}}^\ast(1)\right) \geq \cdots  \\
&\geq g_i^M\left({\U_i^{(M)}}^\ast(1), {\bm \xi_{-i}^{(M)}}^\ast(1)\right) \geq \cdots \\
&\geq g_i^1\left({\U_i^{(1)}}^\ast(k),{\bm \xi_{-i}^{(1)}}^\ast(k)\right) \geq \cdots \\
&\geq g_i^M\left({\U_i^{(M)}}^\ast(k), {\bm \xi_{-i}^{(M)}}^\ast(k)\right) = \gamma_i'.
\end{aligned}
\end{equation*}

Since $f_i$ is bounded below, the sequence $f_i$ is monotonically non-increasing. By the Monotone Convergence Theorem, it converges to a unique limit point $\gamma_i= \gamma_i'\triangleq \bar{\gamma_i}$. Formally, as $k \to \infty$, $\lim_{k \to \infty} g_i^j \left( \mathbf{U}_i^{(j)\ast}(k), \bm{\xi}_{-i}^{(j)\ast}(k) \right) = \gamma_i$.

The overall algorithm $\Omega$ can be decomposed into sub-algorithms $\Omega = \Omega_1 \circ \cdots \circ \Omega_M $, where each $\Omega_j$ corresponds to solving $g_i^j$. $\Omega$ is a closed algorithm since all the updates are performed by continuous functions. All sub-algorithms
$ g_i^1\left({\U_i^{(1)}}^\ast(k),{\bm\xi_{-i}^{(1)}}^\ast(k)\right)$ decrease the value of $ f_i$, therefore it is clear that \( \Omega \) is monotonic with respect to $f_i$. By the properties of monotonic algorithms combined with the closed nature ensured by continuous updates, $\Omega$
converges.

\section{Numerical experiments}
To verify the effectiveness of the proposed optimization model for linear binary classification problems, we selected several types of datasets from kaggle public databases and conducted a series of numerical experiments. These experiments allowed us to thoroughly investigate the model's performance in terms of classification accuracy and computational efficiency.

All numerical experiments were implemented in MATLAB 9.0 on a personal computer with AMD Ryzen 7 4800H CPU 2.90GHz and 16 GB random-access memory(RAM). For simplicity, we use TC and FC to denote the numbers of true classification and false classification, respectively, and use ACCU to denote the accuracy of classification, i.e., ACCU=${\rm\frac{TC}{TC+FC}}$.
The numerical analysis will focus on the classifying accuracy.

In our numerical experiments, we will choose 5 widely-used classifiers, i.e., the classifiers developed respectively by \cite{cvv95}, \cite{sz11}, \cite{dx07},  \cite{zgw09} and \cite{sz16} to make a numerical comparison with the classifier established in this paper. To ensure the reliability of the experimental results, we employed 10-fold cross-validation and repeated the experiment ten times, with the final results averaged to minimize random variance.

The following table lists some information of these models when doing the numerical experiments, where $\Theta = \{2^{-5},2^{-3}, ...,2^{7}\}$.

\begin{table}[!ht]
\centering
	\caption{classiciation models, parameters and Reference }
    \scalebox{1.0}{
	\begin{tabular}{l|lccc}\hline
		Model &   Method   & Number of Parameter& Parameter Set   & Reference\\\hline
		HSVM & {\rm LIBSVM}      &1 &$\Theta$  & \cite{cvv95} \\
		TWSVM & {\rm Algorithm~1}    &2  &$\Theta$ &\cite{sz11} \\
		STL  & {\rm Algorithm~1}  &1 &$\Theta$  &\cite{dx07} \\
		TWSTM & {\rm Algorithm~1} &2 &$ \Theta$  & \cite{zgw09} \\
		TBSTM & {\rm Algorithm~1}  &4&$\Theta$  &\cite{sz16} \\
		LDM-NPSTM & {\rm Algorithm~1}  &8 &$\Theta$  &This paper\\
		\hline
	\end{tabular}}
	\label{tab01}
\end{table}

\subsection{Numerical experiments on skin lesions}
In this experiment, we selected four types of skin cancer \footnote{https://www.kaggle.com/datasets/farjanakabirsamanta/skin-cancer-dataset} with significant pathological differences for evaluation, namely actinic keratosis, dermatofibroma, squamous cell carcinoma, and vascular lesion. All sample images were converted to grayscale before the experiments.
 Related information of these skin cancer datasets is described in Table \ref{tab02} where Size denotes the number of samples, Features denotes the number of features, Ratio denotes this category's proportion of the total.
\begin{table}[!ht]
	\centering
	\caption{The detail of skin cancer datasets}
		\begin{tabular}{llccccccc}\hline
			Data Set & Data set & Features & Size & Ratio  \\ \hline
			\texttt{Acti} & Actinic keratosis & 600$\times$450 & 114 & 0.215  \\
			\texttt{Derm} & Dermatofibroma & 600$\times$450 & 95 & 0.180  \\
			\texttt{Squa} & Squamous cell carcinoma &600$\times$450  & 181 & 0.343  \\
			\texttt{Vasc} & Vascular lesion &600$\times$450 & 139 & 0.263  \\\hline
	\end{tabular}
	\label{tab02}
\end{table}
To further assess classification performance between different skin cancer types, all possible pairwise combinations of the four skin cancer categories are constructed, and binary classification tasks are conducted accordingly. The numerical results for the datasets are given in Table \ref{tab03}, where {\bf W}, {\bf T} and {\bf L} denote the number of wins, ties and losses against another classifier, respectively.

\begin{table}[!ht]
    \centering
    	\caption{Numerical Comparisons of 6 Solvers for Skin Cancer Datasets}
    	\renewcommand{\arraystretch}{1}
	\resizebox{1.0\columnwidth}{!}{
    \begin{tabular}{c|ccccccccc}
    \hline
        Data Sets &  HSVM & TWSVM & STL & TWSTM & TBSTM & LDM-NPSTM \\ \hline
        (\texttt{Acti, Derm}) & 53.69$\pm$1.85  & 53.69$\pm$1.63  & \textbf{65.72}$\pm$4.38  & 63.50$\pm$1.73  & 62.79$\pm$4.59  & 63.79$\pm$2.40  \\
        &3.7486  & 19.1810  & 10.2615  & 5.0595  & 10.4248  & 4.6033  \\
        (\texttt{Acti, Squa}) & 62.70$\pm$0.89  & 62.70$\pm$0.94  & 71.21$\pm$2.69  & 71.62$\pm$0.66  & 67.53$\pm$2.99  & \textbf{88.35}$\pm$1.79  \\
        &5.4850  & 25.8718  & 10.0563  & 5.0371  & 10.3809  & 3.6986  \\
        (\texttt{Acti, Vasc}) & 43.38$\pm$7.84  & 43.15$\pm$1.20  & 70.79$\pm$2.87  & 79.85$\pm$2.41  & 90.84$\pm$1.47  &\textbf{92.24}$\pm$1.95  \\
        &4.6091  & 22.3561  & 10.3062  & 4.7775  & 10.5695  & 4.6224  \\
        (\texttt{Derm, Squa}) & 41.80$\pm$1.50  & 56.28$\pm$6.27  & 67.31$\pm$0.23  & 67.31$\pm$0.23  & 67.31$\pm$4.02  & \textbf{74.54}$\pm$0.28  \\
        &5.0635  & 19.1270  & 10.2640  & 4.6862  & 10.4803  & 5.7556  \\
        (\texttt{Derm, Vasc}) & 59.51$\pm$4.23  & 55.34$\pm$1.67  & 65.82$\pm$5.31  & 55.41$\pm$1.18  & 76.45$\pm$1.11  & \textbf{78.94}$\pm$3.71  \\
        &4.3154  & 37.4657  & 10.3924  & 4.7324  & 10.3661  & 4.6808  \\
        (\texttt{Squ, Vasc}) & 53.13$\pm$8.84  & 43.75$\pm$4.42  & 81.34$\pm$5.21  & 68.67$\pm$0.29  & 67.85$\pm$4.17  & \textbf{83.98}$\pm$0.41  \\
        &5.6678  & 20.3616  & 10.3448  & 4.8153  & 10.3730  & 3.7172 \\ \hline
      {\bf W/T/L} &6/0/0    &6/0/0      &5/0/1    &6/0/0     &6/0/0 \\ \hline
    \end{tabular}}
\label{tab03}
\end{table}
Table \ref{tab03} indicates that LDM-NPSTM generally achieves higher classification accuracy across the datasets compared to HSVM, TWSVM, STL, TWSTM, and TBSTM. It also tends to exhibit lower time costs, smaller standard deviations, and fewer losses in the W/T/L statistics, suggesting better overall effectiveness, efficiency, and robustness.

\subsection{Numerical experiments on on color image datasets}
In this subsection, we evaluate the performance of the proposed algorithm on color image datasets, focusing on classification accuracy and computational efficiency. A dataset of horse breed \footnote{https://www.kaggle.com/datasets/olgabelitskaya/horse-breeds} images is employed, where each image is represented as a third-order tensor $\mathcal{ X}_i \in \mathbb{R}^{n_1\times n_2\times 3}$, with the last dimension corresponding to the RGB channels. Related information of the dataset is described in Table \ref{tab04}.
\begin{table}[!ht]
	\centering
	\caption{The detail of horse speeds datasets}
		\begin{tabular}{llccccccc}\hline
			Dataset & Data set & Features & Size & Ratio  \\ \hline
			\texttt{Akha} & Akhal-Teke & 256$\times$256$\times$3 & 123 & 0.184  \\
			\texttt{Appa} & Appaloosa & 256$\times$256$\times$3  & 105 & 0.157  \\
			\texttt{Orlo} & Orlov Trotter & 256$\times$256$\times$3  & 107 & 0.160  \\
			\texttt{Vlad} & Vladimir Heavy  & 256$\times$256$\times$3  & 37 & 0.055  \\
			\texttt{Perc} &Percheron & 256$\times$256$\times$3  & 56 & 0.0836  \\
			\texttt{Arab} &Arabian & 256$\times$256$\times$3  & 122 & 0.1821  \\
			\texttt{Frie} & Friesian & 256$\times$256$\times$3  & 120 & 0.1791 \\\hline
	\end{tabular}
	\label{tab04}
\end{table}

The experimental results are reported in Table \ref{tab05}. As observed from the table, in the RGB image classification tasks, the LDM-NPSTM algorithm generally exhibits better overall performance compared to the other methods. Across most datasets, it achieves relatively higher average classification accuracies and lower time costs. In particular, according to the W/T/L statistics, LDM-NPSTM records the fewest losses among all methods, further indicating its effectiveness and advantages in image classification tasks.
\begin{table}[!ht]
    \centering
    \caption{Numerical Comparisons of 6 Solvers for Horse Speeds Datasets}
    	\renewcommand{\arraystretch}{1}
	\resizebox{1.0\columnwidth}{!}{
    \begin{tabular}{c|ccccccccc}
    \hline
      Dataset &HSVM & TWSVM & STL & TWSTM & TBSTM & LDM-NPSTM \\ \hline
        (\texttt{Akha, Appa}) & 55.63$\pm$4.89 &55.73$\pm$11.18 &54.02$\pm$3.49 &59.09$\pm$9.55 &55.53$\pm$4.50 &\textbf{59.78}$\pm$5.47  \\
        &3.0059 &25.4503 &3.7839 &11.9574 &11.2630 &4.0839  \\
        (\texttt{Akha, Orlo}) & 50.00$\pm$3.07 &\textbf{65.22}$\pm$13.95 &56.52$\pm$6.74 &57.25$\pm$6.97 &60.14$\pm$6.97 &63.04$\pm$5.61  \\
        &3.1193 &23.9600 &3.9256 &13.0914 &11.3513 &4.3690  \\
        (\texttt{Akha, Vlad}) & 59.38$\pm$13.26 &68.75$\pm$17.68 &73.96$\pm$2.55 &76.04$\pm$4.70 &76.04$\pm$2.55 &\textbf{76.56}$\pm$3.13  \\
        &1.9063 &14.7929 &6.7207 &8.8831 &5.5984 &3.4126  \\
        (\texttt{Akha, Perc}) &40.20$\pm$9.71 &65.36$\pm$17.56 &65.36$\pm$4.75 &67.16$\pm$8.98 &67.32$\pm$1.60 &\textbf{67.65}$\pm$1.96  \\
        &2.3292 &19.5092 &4.8470 &9.9372 &6.5611 &4.3991  \\
        (\texttt{Akha, Arab}) &53.17$\pm$7.31 &46.92$\pm$1.53 &61.08$\pm$6.92 &55.69$\pm$4.45 &\textbf{62.44}$\pm$7.27 &59.00$\pm$2.97  \\
        &3.2237 &17.2410 &3.8087 &14.1637 &13.6802 &6.1761  \\
        (\texttt{Akha, Frie}) &\textbf{81.75}$\pm$8.13 &55.25$\pm$10.25 &63.89$\pm$6.76 &65.92$\pm$9.84 &74.81$\pm$10.88 &74.54$\pm$3.39  \\
        &3.0481 &19.8162 &3.9888 &14.0261 &12.2812 &4.1400  \\
        (\texttt{Appa, Orlo}) &48.81$\pm$1.68 &55.63$\pm$11.33 &58.66$\pm$7.62 &62.55$\pm$8.13 &\textbf{68.72}$\pm$12.34 &68.29$\pm$2.88  \\
        &2.8456 &23.0162 &3.3459 &11.5770 &9.9158 &6.5099  \\
        (\texttt{Appa, Vlad}) &61.43$\pm$26.26 &72.14$\pm$11.11 &76.98$\pm$10.60 &73.33$\pm$4.62 &73.25$\pm$2.77 &\textbf{79.40}$\pm$5.09  \\
        &1.8655 &21.1660 &4.8944 &5.0748 &3.3881 &1.6993  \\
        (\texttt{Appa, Perc}) &60.48$\pm$5.98 &57.54$\pm$1.82 &56.74$\pm$6.40 &63.85$\pm$7.60 &\textbf{65.99}$\pm$3.13 &63.14$\pm$4.12  \\
        &2.1626 &21.0468 &3.0650 &5.9603 &4.1008 &2.2039  \\
        (\texttt{Appa, Arab}) & 55.73$\pm$11.18 &\textbf{66.50}$\pm$10.48 &55.47$\pm$4.28 &59.85$\pm$5.78 &62.81$\pm$7.44 &62.70$\pm$5.10  \\
        &2.9869 &21.8005 &2.9025 &8.9155 &7.4888 &6.4827  \\
        (\texttt{Appa, Frie}) & \textbf{88.74}$\pm$9.78 &70.91$\pm$12.86 &69.30$\pm$5.78 &84.75$\pm$8.97 &79.58$\pm$8.94 &88.04$\pm$4.16  \\
        &2.6760 &10.8180 &3.3446 &8.9030 &6.4633 &3.1007  \\
        (\texttt{Orlo, Vlad}) & 72.38$\pm$1.35 &\textbf{82.62}$\pm$5.72 &69.29$\pm$5.74 &79.52$\pm$10.48 &73.89$\pm$2.42 &74.64$\pm$2.62  \\
        &1.8443 &17.6141 &3.9880 &4.8473 &3.4914 &1.5625  \\
        (\texttt{Orlo, Perc}) & 57.35$\pm$10.40 &66.73$\pm$2.86 &\textbf{67.77}$\pm$7.62 &62.44$\pm$12.22 &64.64$\pm$2.28 &65.63$\pm$3.47  \\
        &2.2092 &13.9050 &3.0108 &6.0591 &4.3135 &2.8230  \\
        (\texttt{Orlo, Arab}) & 44.47$\pm$1.40 &53.36$\pm$1.68 &57.61$\pm$5.79 &62.88$\pm$8.24 &56.23$\pm$5.26 &\textbf{66.52}$\pm$4.38  \\
        &3.0544 &20.7566 &2.9684 &9.0866 &7.5804 &6.4946  \\
        (\texttt{Orlo, Frie}) & 66.90$\pm$14.67 &55.53$\pm$1.40 &64.23$\pm$6.96 &62.81$\pm$8.39 &\textbf{76.02}$\pm$9.70 &69.17$\pm$6.56  \\
        &2.8396 &18.5657 &2.9819 &8.9926 &7.3749 &5.1969  \\
        (\texttt{Vlad, Perc}) & 57.78$\pm$3.14 &68.89$\pm$12.57 &68.33$\pm$6.58 &68.33$\pm$9.63 &71.85$\pm$5.69 &\textbf{74.17}$\pm$10.67  \\
        &1.3725 &17.2874 &1.8888 &3.4224 &2.0750 &1.1853  \\
        (\texttt{Vlad, Arab}) & 74.17$\pm$1.18 &\textbf{83.96}$\pm$3.83 &78.93$\pm$5.11 &75.90$\pm$9.83 &83.33$\pm$10.21 &82.38$\pm$2.98  \\
        &2.3359 &18.1972 &5.1279 &6.1607 &1.7887 &2.1546  \\
        (\texttt{Vlad, Frie}) & 57.92$\pm$6.48 &74.58$\pm$17.09 &74.86$\pm$7.65 &73.89$\pm$11.99 &71.74$\pm$10.02 &\textbf{77.81}$\pm$3.29  \\
        &2.3546 &17.9876 &4.7237 &6.1368 &2.3467 &1.7795  \\
        (\texttt{Perc, Arab}) &51.14$\pm$14.10 &\textbf{74.35}$\pm$3.00 &68.25$\pm$2.50 &66.23$\pm$9.96 &67.32$\pm$3.86 &73.20$\pm$3.15  \\
        &2.4354 &28.1746 &3.6957 &6.7426 &4.8790 &3.1076  \\
        (\texttt{Perc, Frie}) &71.41$\pm$1.16 &\textbf{77.29}$\pm$7.16 &67.37$\pm$6.09 &67.37$\pm$4.98 &69.12$\pm$9.28 &71.81$\pm$4.61  \\
        &2.4466 &23.4560 &3.5893 &6.7311 &5.5610 &3.0582  \\
        (\texttt{Arab, Frie}) &\textbf{75.42}$\pm$6.48 &63.33$\pm$4.71 &62.42$\pm$6.90 &60.36$\pm$8.78 &72.61$\pm$6.69 &69.50$\pm$6.35  \\
        &3.0436 &18.8568 &2.9633 &9.7408 &8.0380 &2.0672 \\ \hline
        W/T/L  &18/0/3 &14/0/7  &19/0/2 &19/0/2 &13/0/8    \\ \hline
    \end{tabular}}
    \label{tab05}
\end{table}

Next, the numerical results are statistically compared using the Friedman test and the Nemenyi post hoc test \cite{dj06}. The Friedman test is first applied to the ACCU values of different classifiers to evaluate whether significant differences exist. Let $ R_j $ denote the rank of the $ j $-th algorithm over $ N $ datasets. Under the null hypothesis that all algorithms perform equivalently, their average ranks $ R_j $ should be equal. The Friedman statistic is computed as

$$
\chi_F^2 = \frac{12N}{K(K + 1)} \left( \sum_{j=1}^K R_j^2 - \frac{K(K + 1)^2}{4} \right).
$$

If the null hypothesis is rejected, the Nemenyi post hoc test is conducted to identify specific differences. Two classifiers are considered significantly different if the absolute difference in their average ranks exceeds the critical difference (CD), given by
$$
{\rm CD} = q_{\alpha} \sqrt{ \frac{K(K+1)}{6N} },
$$
where $K$ and $N$ denote the number of algorithms and datasets, respectively, and $q_{\alpha}$ is the critical value of the Studentized range statistic based on a given signicance level and this value can be obtained from the threshold table of $q_{\alpha}$ (see Table 5 in \cite{dj06}).

\begin{figure*}[!ht]
	\includegraphics[width=0.45\textwidth]{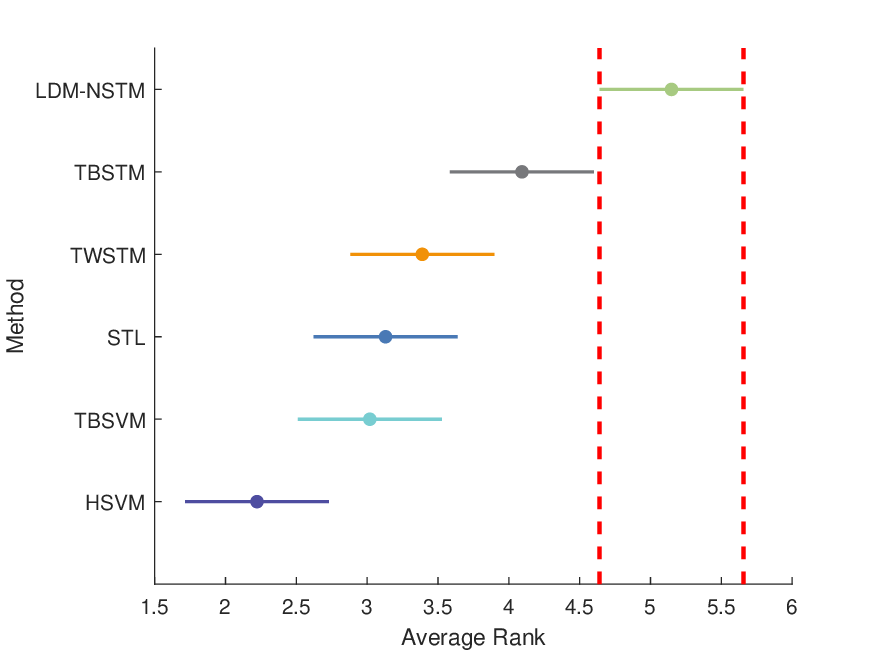}
	\quad\includegraphics[width=0.45\textwidth]{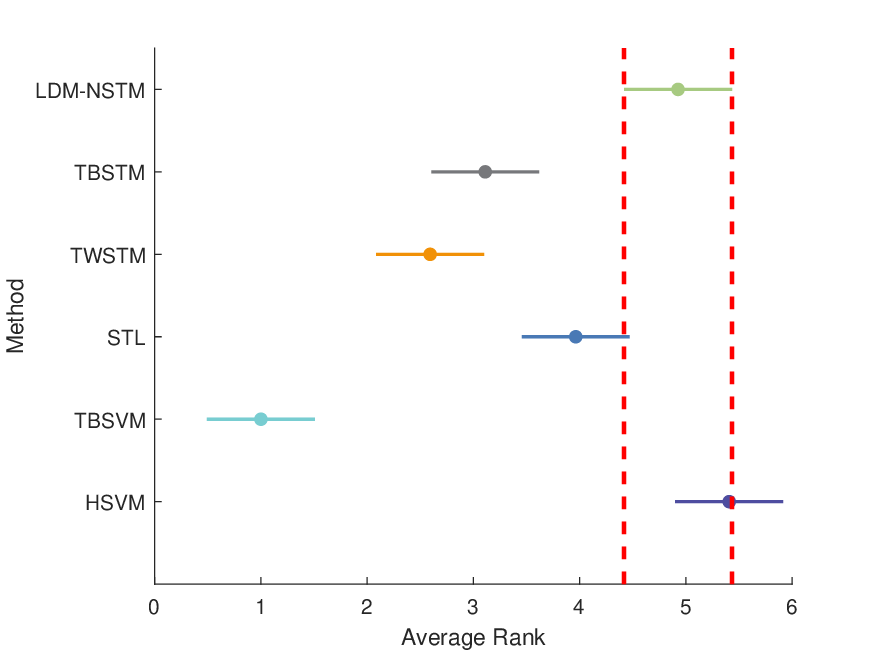}
	\caption{Nemenyi posthoc test with ACCU (left) and CPU Time (right)}
	\label{f501}
\end{figure*}

In the numerical experiments of this subsection, 6 classifiers are evaluated across 21 datasets, and their rankings based on ACCU are reported in Table \ref{tab06}. Using the average ranks $ R_j $ provided in the last row of the table, the Friedman test yields a $P$-value of $ 1.78\times10^{-7} $. Given the commonly adopted significance level $\alpha = 0.05$. Surely, the value is signicantly less than the signicant level, this result indicates significant differences among the classifiers. Substituting $K=6$, $N=21$, and $\alpha=0.05$ into the formula for critical difference (CD) results in ${\rm CD} = 0.5081$. Based on this, the Nemenyi post hoc comparison results are shown on the left side of Fig. \ref{f501}. It is evident that the absolute differences in the mean ranks between our model and other models, indicating statistically significant differences in ACCU between LDM-NSVM and most other models.

Similarly, Table \ref{tab07} presents the rankings of the CPU times for the 6 classifiers across the 27 datasets. In this case, the Friedman test yields a $ P $-value of $ 3.03\times10^{-17} $, with the same CD value of $ 0.5081 $. The corresponding Nemenyi post hoc results are shown on the right side of Fig. \ref{f501}. It can be observed that the absolute differences in mean ranks between LDM-NSVM and all other models exceed the CD value except HSVM, demonstrating significant statistical differences in CPU time between LDM-NSVM  and most other models.

\begin{table}
    \centering
    \caption{Average rank of 6 Solvers for ACCU}
	\renewcommand{\arraystretch}{0.85}
    \resizebox{0.6\columnwidth}{!}{
    \begin{tabular}{c|cccccc}
     \hline
      Dataset& HSVM& TWSVM& STL& TWSTM& TBSTM& LDM \_NPSTM \\\hline
         (\texttt{Acti, Derm})& 1.5 & 1.5 & 6 & 4 & 3 & 5 \\
         (\texttt{Acti, Squa})& 1.5 & 1.5 & 4 & 5 & 3 & 6 \\
         (\texttt{Acti, Vasc})& 2 & 1 & 3 & 4 & 5 & 6 \\
         (\texttt{Derm, Squa})& 1 & 2 & 4 & 4 & 4 & 6 \\
         (\texttt{Derm, Vasc})& 3 & 1 & 4 & 2 & 5 & 6 \\
         (\texttt{Squ, Vasc}) & 2 & 1 & 5 & 4 & 3 & 6 \\
         (\texttt{Akha, Appa})& 3 & 4 & 1 & 5 & 2 & 6 \\
         (\texttt{Akha, Orlo})& 1 & 6 & 2 & 3 & 4 & 5 \\
         (\texttt{Akha, Vlad})& 1 & 2 & 3 & 4.5 & 4.5 & 6 \\
         (\texttt{Akha, Perc})& 1 & 2.5 & 2.5 & 4 & 5 & 6 \\
         (\texttt{Akha, Arab})& 2 & 1 & 5 & 3 & 6 & 4 \\
         (\texttt{Akha, Frie})& 6 & 1 & 2 & 3 & 5 & 4 \\
         (\texttt{Appa, Orlo})& 1 & 2 & 3 & 4 & 6 & 5 \\
         (\texttt{Appa, Vlad})& 1 & 2 & 5 & 4 & 3 & 6 \\
         (\texttt{Appa, Perc})& 3 & 2 & 1 & 5 & 6 & 4 \\
         (\texttt{Appa, Arab})& 2 & 6 & 1 & 3 & 5 & 4 \\
         (\texttt{Appa, Frie})& 6 & 2 & 1 & 4 & 3 & 5 \\
         (\texttt{Orlo, Vlad})& 2 & 6 & 1 & 5 & 3 & 4 \\
         (\texttt{Orlo, Perc})& 1 & 5 & 6 & 2 & 3 & 4 \\
         (\texttt{Orlo, Arab})& 1 & 2 & 4 & 5 & 3 & 6 \\
         (\texttt{Orlo, Frie})& 4 & 1 & 3 & 2 & 6 & 5 \\
         (\texttt{Vlad, Perc})& 1 & 4 & 2.5 & 2.5 & 5 & 6 \\
         (\texttt{Vlad, Arab})& 1 & 6 & 3 & 2 & 5 & 4 \\
         (\texttt{Vlad, Frie})& 1 & 4 & 5 & 3 & 2 & 6 \\
         (\texttt{Perc, Arab})& 1 & 6 & 4 & 2 & 3 & 5 \\
         (\texttt{Perc, Frie})& 4 & 6 & 1.5 & 1.5 & 3 & 5 \\
         (\texttt{Arab, Frie})& 6 & 3 & 2 & 1 & 5 & 4 \\\hline
        $R_j$ & 2.2222  & 3.0185  & 3.1296  & 3.3889  & 4.0926  & 5.1481 \\ \hline
    \end{tabular}}
      \label{tab06}
\end{table}

\begin{table}
    \centering
    \caption{Average rank of 6 Solvers for CPU Time}
    \renewcommand{\arraystretch}{0.85}
	\resizebox{0.6\columnwidth}{!}{
    \begin{tabular}{c|ccccccc}\hline
        Dataset& HSVM&TWSVM&STL&TWSTM&TBSTM&LDM\_NPSTM \\\hline
        (\texttt{Acti, Derm}) & 6 & 1 & 3 & 4 & 2 & 5 \\
        (\texttt{Acti, Squa}) & 4 & 1 & 3 & 5 & 2 & 6 \\
        (\texttt{Acti, Vasc}) & 6 & 1 & 3 & 4 & 2 & 5 \\
        (\texttt{Derm, Squa}) & 5 & 1 & 3 & 6 & 2 & 4 \\
        (\texttt{Derm, Vasc}) & 6 & 1 & 2 & 4 & 3 & 5 \\
        (\texttt{Squ, Vasc})  & 4 & 1 & 3 & 5 & 2 & 6 \\
        (\texttt{Akha, Appa}) & 6 & 1 & 5 & 2 & 3 & 4 \\
        (\texttt{Akha, Orlo}) & 6 & 1 & 5 & 2 & 3 & 4 \\
        (\texttt{Akha, Vlad}) & 6 & 1 & 3 & 2 & 4 & 5 \\
        (\texttt{Akha, Perc}) & 6 & 1 & 4 & 2 & 3 & 5 \\
        (\texttt{Akha, Arab}) & 6 & 1 & 5 & 2 & 3 & 4 \\
        (\texttt{Akha, Frie}) & 6 & 1 & 5 & 2 & 3 & 4 \\
        (\texttt{Appa, Orlo}) & 6 & 1 & 5 & 2 & 3 & 4 \\
        (\texttt{Appa, Vlad}) & 5 & 1 & 3 & 2 & 4 & 6 \\
        (\texttt{Appa, Perc}) & 6 & 1 & 4 & 2 & 3 & 5 \\
        (\texttt{Appa, Arab}) & 5 & 1 & 6 & 2 & 3 & 4 \\
        (\texttt{Appa, Frie}) & 6 & 1 & 4 & 2 & 3 & 5 \\
        (\texttt{Orlo, Vlad}) & 5 & 1 & 3 & 2 & 4 & 6 \\
        (\texttt{Orlo, Perc}) & 6 & 1 & 4 & 2 & 3 & 5 \\
        (\texttt{Orlo, Arab}) & 5 & 1 & 6 & 2 & 3 & 4 \\
        (\texttt{Orlo, Frie}) & 6 & 1 & 5 & 2 & 3 & 4 \\
        (\texttt{Vlad, Perc}) & 5 & 1 & 4 & 2 & 3 & 6 \\
        (\texttt{Vlad, Arab}) & 4 & 1 & 3 & 2 & 6 & 5 \\
        (\texttt{Vlad, Frie}) & 4 & 1 & 3 & 2 & 5 & 6 \\
        (\texttt{Perc, Arab}) & 6 & 1 & 4 & 2 & 3 & 5 \\
        (\texttt{Perc, Frie}) & 6 & 1 & 4 & 2 & 3 & 5 \\
        (\texttt{Arab, Frie}) & 4 & 1 & 5 & 2 & 3 & 6 \\\hline
         $R_j$& 5.4074  & 1.0000  & 3.9630  & 2.5926  & 3.1111 & 4.9259 \\ \hline
    \end{tabular}}
      \label{tab07}
\end{table}

\subsection{Numerical experiments on  small sample classification problems}
In this experiment, we used musical instruments dataset \footnote{\footnotesize https://www.kaggle.com/datasets/nikolasgegenava/music-instruments} to validate the effectiveness of our model on small sample classification problems. Musical instruments dataset is an image dataset, which includes 30 object categories of different musical instruments.
Note that the number of images in each category differs a lot, about 70 to 200 images per category. However, the size of the images in each category is equal to 200$\times$300$\times$3. Then, we selected six types of musical instruments, namely \texttt{drums},\texttt{harp} \texttt{piano},\texttt{sax}, \texttt{sitar} and \texttt{violin}, for binary classification experiments, with each category consisting of 60 color images.

Table \ref{tab08} respectively lists the classification results of the five comparison methods and the proposed model. We observe that our model achieves almost similar or better performance than the baseline in all class pairs, which indicates that our model has potential benefits when dealing with small sample classification problems.

\begin{table}[!ht]
    \centering
    \caption{Numerical Comparisons of 6 Solvers for musical instruments Datasets}
    \renewcommand{\arraystretch}{0.9}
    \resizebox{0.85\columnwidth}{!}{
    \begin{tabular}{lccccccccc}
    \hline
         Dataset& HSVM&TWSVM&STL&TWSTM&TBSTM&LDM\_NPSTM \\\hline
        (\texttt{drums,harp}) & 56.25$\pm$8.84&74.66$\pm$5.79&68.33$\pm$9.13&66.25$\pm$13.69&72.92$\pm$ 14.73& {\bf77.08}$\pm$8.84 \\
       &0.6985&19.1593&4.9143&4.9499&0.1528&0.1240    \\
        (\texttt{drums,piano}) & 53.52$\pm$15.98&63.73$\pm$5.98&55.00$\pm$7.45&68.75$\pm$10.87&68.75$\pm$2.95& {\bf72.92}$\pm$2.95 \\
       &0.6466&24.2987&4.8740&5.1063&0.2381&0.2265    \\
        (\texttt{drums,sax}) & 68.75$\pm$8.84&78.54$\pm$5.61&76.67$\pm$14.91&77.50$\pm$13.94&84.58$\pm$ 2.95& {\bf87.42}$\pm$2.95 \\
       &0.7209&28.0308&5.0962&5.1396&0.1667&0.1417    \\
        (\texttt{drums,sitar})&75.00$\pm$8.71&73.18$\pm$6.02&68.33$\pm$9.13&70.00$\pm$17.08&{\bf68.75}$\pm$2.95& 67.55$\pm$4.73 \\
       &0.7480&25.1164&5.1823&4.9459&0.1573&0.1438    \\
        (\texttt{drums,violin}) & 66.67$\pm$11.79&60.84$\pm$6.00&66.67$\pm$8.33&61.25$\pm$9.13&62.50$\pm$5.89& {\bf68.58}$\pm$8.84 \\
       &0.8101&23.2174&5.5205&5.1270&0.1739&0.2381    \\
        (\texttt{harp,piano}) & 60.42$\pm$14.73&58.92$\pm$5.98&61.67$\pm$7.45&57.50$\pm$5.89&64.58$\pm$ 8.84& {\bf 68.75}$\pm$2.95 \\
       &0.7836&23.3995&5.0960&5.1261&0.2314&0.1892    \\
        (\texttt{harp,sax}) & 56.25$\pm$20.62&{\bf73.17}$\pm$6.03&70.00$\pm$9.50&64.50$\pm$12.36&61.50$\pm$ 5.89& 68.92$\pm$11.79 \\
       &0.7595&24.3436&5.7827&5.1208&0.1866&0.2173    \\
        (\texttt{harp,sitar})&45.83$\pm$11.79&53.33$\pm$6.41&{\bf66.67}$\pm$8.33&58.75$\pm$9.50&56.25$\pm$ 8.84& 64.17$\pm$5.89 \\
       &0.7763&26.3310&6.8348&5.1513&0.2218&0.2074    \\
        (\texttt{harp,violin}) & 45.83$\pm$5.89&71.84$\pm$6.09&65.00$\pm$6.97&62.50$\pm$6.97&73.08$\pm$ 2.95& {\bf77.35}$\pm$14.73 \\
       &0.7788&27.8693&8.3216&5.1590&0.1924&0.1884    \\
         (\texttt{piano,sax}) & 41.67$\pm$23.57&{\bf86.67}$\pm$6.11&83.33$\pm$5.89&73.75$\pm$9.50&75.00$\pm$5.89& 72.92$\pm$14.73 \\
       &0.8405&23.8435&9.3242&5.0618&0.1258&0.1101    \\
        (\texttt{piano,sitar})&{\bf77.08}$\pm$2.95&63.33$\pm$6.46&71.67$\pm$7.45&63.75$\pm$11.18&68.75$\pm$14.73& 66.67$\pm$5.89 \\
       &0.7850&24.2790&8.7510&5.5522&0.2431&0.1818    \\
         (\texttt{piano,violin}) & 25.00$\pm$5.89&71.67$\pm$6.11&71.33$\pm$14.91&62.50$\pm$12.64&58.33$\pm$11.79& {\bf72.50}$\pm$11.51 \\
       &0.7725&25.7401&8.5646&5.2260&0.1776&0.1545    \\
        (\texttt{sax,sitar})&66.67$\pm$5.89&65.00$\pm$6.14&68.33$\pm$16.03&{\bf73.75}$\pm$6.97&70.83$\pm$ 17.68& 68.75$\pm$14.73 \\
       &0.8401&28.4360&8.3660&5.2314&0.2346&0.2231    \\
        (\texttt{sax,violin}) & 52.08$\pm$2.95&63.18$\pm$6.20&63.74$\pm$11.18&61.25$\pm$5.89&56.25$\pm$ 8.84& {\bf64.58}$\pm$14.73 \\
       &0.7983&24.4178&8.4763&5.1999&0.2019&0.1759    \\
         (\texttt{sitar,violin}) & 60.83$\pm$5.89&61.67$\pm$6.11&60.00$\pm$10.87&58.75$\pm$13.94&58.33$\pm$5.11& {\bf65.42}$\pm$2.95 \\
       &0.8308&26.8856&8.5614&5.4858&0.2359&0.1970    \\\hline
          W/T/L& 13/0/2&12/0/3 &10/0/5 &12/0/3 &11/0/4    \\ \hline
    \end{tabular}}
    \label{tab08}
\end{table}

\section{Conclusion}
In this paper, we proposed a novel LDM-NPSTM for binary classification tasks, which integrates both distributional and sampling information of training data within a tensor-based learning framework. Specifically, the separating hyperplane parameters in LDM-NPSTM form a tensorplane decomposed into a sum of rank-one tensors via CP decomposition, enabling efficient exploitation of multiway structural information. The corresponding optimization problems were solved iteratively using a  alternating projection strategy, where each subproblem along a specific tensor mode was formulated as a standard SVM-type convex optimization problem. The efficiency of the proposed method
was illustrated on the problems of skin lesions datasets, color image datasets, and small sample classification.



\begin{thebibliography}{plain}
\bibitem{an10} Amayri, O., Nizar, B.: A study of spam filtering using support vector machines. Artificial Intelligence Review. 34, 73-108 (2010).
\bibitem{bc98} Burges, C, JC.: A tutorial on support vector machines for pattern recognition. Data mining and knowledge discovery. 2(2), 121-167 (1998).
\bibitem{by10} Bollegala, D., Yutaka, M., Mitsuru, I.: A web search engine-based approach to measure semantic similarity between words. IEEE Transactions on knowledge and Data Engineering. 23(7), 977-990 (2010).
\bibitem{cf03} Cao, L., Francis, E.: Support vector machine with adaptive parameters in financial time series forecasting. IEEE Transactions on neural networks. 14(6), 1506-1518 (2003).
\bibitem{cc22} Chen C, et al. Kernelized support tensor train machines. Pattern Recognition. 122, 108337 (2022).
\bibitem{cvv95} Cortes, C., Vladimir V.: Support-vector networks. Machine learning. 20, 273-297 (1995).
\bibitem{dl97} De, L., Lieven.: Signal processing based on multilinear algebra. Leuven: Katholieke Universiteit Leuven. (1997).
\bibitem{dj06} Dem$\check{s}$ar J.: Statistical comparisons of classifiers over multiple data sets. Journal of Machine Learning Research. 7(1), 1-30 (2006).
\bibitem{ds08} Dhanjal, C., Steve, R., John S.: Efficient sparse kernel feature extraction based on partial least squares. IEEE Transactions on Pattern Analysis and Machine Intelligence. 31(8), 1347-1361 (2008).
\bibitem{kr97} Etemad, K., Chellappa, R.: Discriminant analysis for recognition of human face images.Journal of the Optical Society of America A, 14(8), 1724-1733 (1997).
\bibitem{gz13} Gao, W., Zhou, Z. H.: On the doubt about margin explanation of boosting. Artificial Intelligence. 203, 1-18 (2013).
\bibitem{gd99} Gavrila D M. The visual analysis of human movement: A survey. Computer vision and image understanding. 73(1), 82-98 (1999).
\bibitem{rl07} Green, R. D., Guan, L.: Quantifying and recognizing human movement patterns from monocular video images-part II: applications to biometrics. IEEE Transactions on Circuits and Systems for Video Technology. 14(2), 191-198 (2004).
\bibitem{hz20} He Z, et al. Support tensor machine with dynamic penalty factors and its application to the fault diagnosis of rotating machinery with unbalanced data. Mechanical systems and signal processing. 141, 106441 (2020).
\bibitem{il08} Isa, D., Lee, L. H., Kallimani, V. P., Rajkumar, R. Text document preprocessing with the Bayes formula for classification using the support vector machine. IEEE Transactions on Knowledge and Data engineering.  20(9), 1264-1272 (2008).
\bibitem{krs07} Khemchandani, R., Chandra, S.: Twin support vector machines for pattern classification. IEEE Transactions on pattern analysis and machine intelligence. 29(5), 905-910 (2007).
\bibitem{kk13} Khemchandani, R., Karpatne, A., Chandra, S.: Proximal support tensor machines, International Journal of Machine Learning and Cybernetics, 4, 703-712 (2013).
\bibitem{km01} Kim M, et al. Moving object segmentation in video sequences by user interaction and automatic object tracking. Image and Vision Computing. 19(5), 245-260 (2001).
\bibitem{kt09} Kolda, T. G., Bader, B. W.: Tensor decompositions and applications. SIAM review 51 (3), 455-500 (2009).
\bibitem{kt06} Kolda, T. G. Multilinear operators for higher-order decompositions (No. SAND2006-2081). Sandia National Laboratories (SNL), Albuquerque, NM, and Livermore, CA (United States) (2006).
\bibitem{lc08} Li, Y., Guan, C.: Joint feature re-extraction and classification using an iterative semi-supervised support vector machine algorithm. Machine Learning. 71, 33-53 (2008).
\bibitem{hk08} Lu, H., Plataniotis, K. N., Venetsanopoulos, A. N.: multilinear principal component analysis of tensor objects. 18th International Conference on Pattern Recognition. 2, 776-779 (2006).
\bibitem{lp09} Lu, H., Plataniotis, K. N., Venetsanopoulos, A. N.: A taxonomy of emerging multilinear discriminant analysis solutions for biometric signal recognition. Biometrics: Theory, Methods, and Applications. 21-45 (2009).
\bibitem{mw06} Mangasarian, O. L., Wild, E. W.: Multisurface proximal support vector classification via generalized eigenvalues. IEEE transactions on pattern analysis and machine intelligence. 28(1), 69-74 (2006).
\bibitem{nl10} Narwaria, M., Lin, W.: Objective image quality assessment based on support vector regression. IEEE Transactions on Neural Networks. 21(3), 515-519 (2010).
\bibitem{pk01} Plataniotis K N. Color image processing and applications. Measurement Science and Technology. 12(2), 222-222 (2001).
\bibitem{ra11} Rahman, M. M., Antani, S. K., Thoma, G. R. A learning-based similarity fusion and filtering approach for biomedical image retrieval using SVM classification and relevance feedback. IEEE Transactions on information technology in biomedicine. 15(4), 640-646 (2011).
\bibitem{nb09} Renard, N., Bourennane, S.: Dimensionality reduction based on tensor modeling for classification methods. IEEE Transactions on Geoscience and Remote Sensing. 47(4), 1123-1131 (2009).
\bibitem{sj10} Sahbi, H., Audibert, J. Y., Keriven, R.: Context-dependent kernels for object classification. IEEE transactions on pattern analysis and machine intelligence. 33(4), 699-708 (2010).
\bibitem{ss96} Sain, S. R.: The nature of statistical learning theory. 409-409 (1996).
\bibitem{vv13} Vapnik, V.: The nature of statistical learning theory. Springer science business media. (2013).
\bibitem{sz11} Shao, Y. H., et al.: Improvements on twin support vector machines. IEEE transactions on neural networks. 22(6), 962-968 (2011).
\bibitem{sk08} Shen, K. Q., et al.: Feature selection via sensitivity analysis of SVM probabilistic outputs. Machine Learning. 70, 1-20 (2008).
\bibitem{sz16} Shi, H., et al.:  Twin bounded support tensor machine for classification. International Journal of Pattern Recognition and Artificial Intelligence, 30(01), 1650002 (2016).
\bibitem{dx05} Tao, D., Li, X., Hu, W., Maybank, S., Wu, X.: Supervised tensor learning. In Fifth IEEE International Conference on Data Mining. IEEE. pages 8-pp (2005).
\bibitem{dx07} Tao, D., Li, X., Wu, X., Hu, W., Maybank, S. J.: Supervised tensor learning. Knowledge and Information Systems. 13(1), 1-42 (2007).
\bibitem{vv99} Vapnik, V. N.: An overview of statistical learning theory. IEEE transactions on neural networks. 10(5), 988-999 (1999).
\bibitem{wa04} Wang H, Ahuja N. Compact representation of multidimensional data using tensor rank-one decomposition. vectors, 1(5), 44-47 (2004).
\bibitem{lm11} Wang, L., et al.: A refined margin analysis for boosting algorithms via equilibrium margin. The Journal of Machine Learning Research. 12, 1835-1863 (2011).
\bibitem{zl23} Zhang L, et al. A novel dual-center-based intuitionistic fuzzy twin bounded large margin distribution machines. IEEE Transactions on Fuzzy Systems. 31(9): 3121-3134 (2023).
\bibitem{zz14} Zhang, T., Zhou, Z. H. Large margin distribution machine. In Proceedings of the 20th ACM SIGKDD international conference on Knowledge discovery and data mining. 313-322 (2014).
\bibitem{zgw09} Zhang, X., Gao, X., Wang, Y. Twin support tensor machines for MCs detection. Journal of Electronics (China). 26(3), 318-325 (2009).

\end{thebibliography}
\end{document}